\documentclass[11pt,a4paper]{article}
\usepackage{amsmath,amsthm,amssymb,mathrsfs,mathtools}
\usepackage[a4paper,margin=1in]{geometry}
\usepackage{cite}
\usepackage[hidelinks]{hyperref}
\newtheorem{theorem}{Theorem}[section]
\newtheorem{lemma}[theorem]{Lemma}
\newtheorem{proposition}[theorem]{Proposition}
\newtheorem{corollary}[theorem]{Corollary}
\theoremstyle{definition}
\newtheorem{definition}[theorem]{Definition}
\theoremstyle{remark}
\newtheorem{remark}[theorem]{Remark}
\newcommand{\R}{\mathbb R}
\newcommand{\Div}{\operatorname{div}}
\DeclareMathOperator{\mtrace}{tr}
\newcommand{\BUC}{\mathrm{BUC}}
\newcommand{\dsym}{D_{\rm sym}}
\newcommand{\ee}{\mathrm e}
\newcommand{\Lop}{\mathcal L}
\newcommand{\bx}{\mathbf e}
\newcommand{\dx}{\,dx}
\numberwithin{equation}{section}
\allowdisplaybreaks[2]
\makeatletter
\renewcommand{\@maketitle}{
  \begin{center}
  {\LARGE\@title\par}
  \ifx\@author\@empty\else
    \vskip 1em{\large\begin{tabular}[t]{c}\@author\end{tabular}\par}
  \fi
  \ifx\@date\@empty\else\vskip .5em{\large\@date\par}\fi
  \end{center}\vskip .6em}
\makeatother
\hypersetup{pdftitle={Global strong solutions to the three-dimensional quantum Navier-Stokes-Poisson system with far-field vacuum}}
\begin{document}
\makeatletter
\renewcommand{\@fnsymbol}[1]{}
\makeatother
\title{Global strong solutions to the three-dimensional quantum
Navier--Stokes--Poisson system with far-field vacuum}
\author{
	\bf\large Xiangdi Huang$^{a}$\thanks{E-mail addresses: xdhuang@amss.ac.cn (X. Huang); zhouhuitao141@163.com (H. Zhou).}, Huitao Zhou$^{a}$\\
	\small a. Institute of Mathematics, Academy of Mathematics and Systems Science,\\
	\small Chinese Academy of Sciences, Beijing 100190, China;\\
}
\date{}
\maketitle
\begin{abstract}
We prove global existence and uniqueness in a Sobolev class for the repulsive
quantum Navier--Stokes--Poisson system on $\R^3$ with zero background charge,
viscosities $\mu(\rho)=\rho$, $\lambda(\rho)=0$, capillarity
$\kappa(\rho)=1/\rho$, and pressure $\rho^\gamma$, $1\le\gamma\le3$.
The initial density $\rho_0\in L^1\cap H^3$ is everywhere positive, although
this class already implies $\rho_0(x)\to0$ as $|x|\to\infty$, so the far-field
vacuum is genuine and no positive background density is assumed.
We assume $u_0\in H^2$, $g_0=\nabla\log\rho_0\in L^\infty$,
and $\nabla g_0\in H^1$. For $1<\gamma<5/4$ we also require
$\rho_0^{\gamma-1}g_0\in L^2$. The data may be arbitrarily large.
On every finite time interval the density remains pointwise comparable to
$\rho_0$. In the isothermal case $u+tg_0\in C([0,T];H^2)$, whereas
$u(t)\notin L^2$ for every $t>0$ and every admissible datum.
The physical and effective energy balances, together with the Fisher
information identity, control $\int_0^T\|\rho(t)\|_\infty^2\,dt$.
This estimate and mass conservation bound the electric field in
$L^3(0,T;L^\infty)$ and close the effective velocity estimate up to
and including the cubic pressure law.\\[4mm]
\end{abstract}

\section{Introduction}\label{sec:introduction}

We consider a charged quantum fluid of finite mass in $\R^3$ whose density
is positive at every point and tends to zero at infinity. Its viscosity
vanishes with the density, while the electric acceleration depends on the
entire charge distribution. We prove that a strong solution with these
properties exists globally for initial data of arbitrary size in the
class of Theorem~\ref{thm:main}. On each finite time interval its density
retains the spatial decay of the initial profile.

Let $\rho$, $u$ and $V$ denote the density, velocity and electrostatic
potential. With dielectric parameter $\ell>0$ and zero background charge,
the equations are
\begin{subequations}\label{eq:model}
\begin{align}
 \rho_t+\Div(\rho u)&=0,\label{eq:mass-pde}\\
 (\rho u)_t+\Div(\rho u\otimes u)+D\rho^\gamma
 &=\Div(2\rho\dsym u)+\Div(\rho D^2\log\rho)-\rho DV,\label{eq:mom}\\
 -\ell^2\Delta V&=\rho,\qquad \lim_{|x|\to\infty}V(t,x)=0.\label{eq:poisson}
\end{align}
\end{subequations}
Here $D$ denotes spatial differentiation, $(Du)_{ij}=\partial_j u_i$,
$\dsym u=(Du+Du^T)/2$, and $(\Div F)_i=\sum_j\partial_j F_{ij}$.
The force $-\rho DV$ is repulsive with the sign convention in
\eqref{eq:poisson}. The viscosity coefficients are $\mu(\rho)=\rho$ and
$\lambda(\rho)=0$. For smooth positive density the capillary stress
associated with $\kappa(\rho)=1/\rho$ gives the Bohm force
\begin{equation*}
 \Div(\rho D^2\log\rho)=2\rho D\left(\frac{\Delta\sqrt\rho}{\sqrt\rho}\right).
\end{equation*}
Such density-gradient stresses belong to Korteweg's capillary theory
\cite[Section 1]{DS}. A kinetic derivation of viscous corrections to the
isothermal quantum Euler equations is given in \cite{BM}.
The initial and far-field conditions are
\begin{equation}\label{eq:initial-far}
 (\rho,u)|_{t=0}=(\rho_0,u_0),\qquad
 (\rho,u,DV)(t,x)\longrightarrow(0,0,0) \quad (|x|\to\infty).
\end{equation}

\subsection*{Main result}
Set
\begin{equation*}
 h=\log\rho,\qquad g=Dh,\qquad v=u+g,\qquad E=DV .
\end{equation*}
All spatial function spaces are over $\R^3$. For a fixed $T>0$, write
$C_tX=C([0,T];X)$ and $L^p_tX=L^p(0,T;X)$.
On each compact space-time cylinder, a continuous positive density
has a positive minimum. This allows us to use $\log\rho$ locally
even when $\rho$ tends to zero at infinity.

\begin{definition}[Positive strong solution]\label{def:positive-strong}
On $[0,T]$, we call $(\rho,u,V)$ a positive strong solution if $\rho$
is continuous and everywhere positive, the distributional derivatives
in \eqref{eq:model} are represented by locally integrable functions on
$(0,T)\times\R^3$, the equations hold almost everywhere there, and the
far-field condition in \eqref{eq:initial-far} holds at each time.
For every $t\in[0,T]$, we take $V(t,\cdot)$ to be the continuous
representative vanishing at infinity that satisfies
$-\ell^2\Delta V(t,\cdot)=\rho(t,\cdot)$ in $\mathcal D'(\R^3)$.
\end{definition}

To describe the regularity of such solutions, we compare $g$ with its
initial profile $g_0=D\log\rho_0$. For isothermal pressure this profile
also enters the time-dependent correction to the velocity.

\begin{theorem}\label{thm:main}
Let $\ell>0$, $1\le\gamma\le3$, and let the continuous representative of
$\rho_0$ be strictly positive on $\R^3$. Assume
\begin{equation}\label{eq:initial}
 \rho_0\in L^1\cap H^3,\qquad u_0\in H^2,\qquad
 g_0=D\log\rho_0\in L^\infty,\qquad Dg_0\in H^1 .
\end{equation}
For $\gamma>1$, write $p_0=\gamma\rho_0^{\gamma-1}g_0$ for the initial
pressure gradient divided by the density. When $1<\gamma<5/4$, assume
in addition that
\begin{equation}\label{eq:initial-pressure}
 p_0\in L^2 .
\end{equation}
Write $v_0=u_0+g_0$ and let $\sigma=1$ for $\gamma=1$ and $\sigma=0$
for $\gamma>1$. Then \eqref{eq:model}--\eqref{eq:initial-far} admits a
global positive strong solution in the sense of
Definition~\ref{def:positive-strong}. For each $T>0$ it has the regularity
\begin{subequations}\label{eq:regularity}
\begin{align}
 &\rho\in C_tH^3\cap L^2_tH^4,\qquad
   \rho_t\in C_tH^1\cap L^2_tH^2,\label{eq:reg-rho}\\
 &u+\sigma tg_0,\quad v-v_0+\sigma tg_0,\quad g-g_0\in C_tH^2,
   \label{eq:reg-shifts}\\
 &Du,\ Dv,\ Dg\in C_tH^1\cap L^2_tH^2,\qquad
   u\in L^2_tH^3\quad\text{if }\gamma>1,\label{eq:reg-spatial}\\
 &u_t+\sigma g_0,\quad v_t+\sigma g_0,\quad g_t
   \in C_tL^2\cap L^2_tH^1,\label{eq:reg-time}\\
 &E\in C_tH^4\cap L^2_tH^5.\label{eq:reg-field}
\end{align}
\end{subequations}
The initial traces are $\rho(0)=\rho_0$ in $H^3$ and
$(u+\sigma tg_0)|_{t=0}=u_0$ in $H^2$. On $[0,T]$ the solution is
unique among positive strong solutions satisfying \eqref{eq:regularity}
and these traces.

Let $M=\int_{\R^3}\rho_0\,dx$. There is a constant $C_T\ge1$ depending
only on $T,\ell,\gamma$, the norms in \eqref{eq:initial}, the fixed
mollifier below and, for $1<\gamma<5/4$, $\|p_0\|_2$, such that
\begin{align}
 &\int_{\R^3}\rho(t,x)\,dx=M,\qquad
   C_T^{-1}\rho_0(x)\le\rho(t,x)\le C_T\rho_0(x)
   \quad (0\le t\le T,\ x\in\R^3),\label{eq:reg-relative}\\
 &\sup_{0\le t\le T}\bigl(\|\rho(t)\|_\infty+\|v(t)\|_\infty+
   \|g(t)\|_\infty+\|E(t)\|_{W^{1,\infty}}\bigr)\le C_T .
   \label{eq:reg-sup}
\end{align}
The potential belongs to $C_tC_0(\R^3)$.
Moreover, $\rho,u,v,g,E$ and $V$ tend to zero as $|x|\to\infty$,
uniformly for $t\in[0,T]$.
\end{theorem}

Lemma~\ref{lem:initial-buc} shows that \eqref{eq:initial} implies
\begin{equation*}
 g_0\in L^6\cap C_0(\R^3)
\end{equation*}
and that \eqref{eq:initial-pressure} follows from \eqref{eq:initial} for
$5/4\le\gamma\le3$. For every $1<\gamma<5/4$ this implication fails.
Moreover, an $L^2$ trace of $v_t$ at zero in the stated spatial
regularity class requires \eqref{eq:initial-pressure}
(Remark~\ref{rem:pressure-threshold}).

At $\gamma=1$, the velocity trace in \eqref{eq:reg-shifts} involves
$u+tg_0$. To determine what this says about $u(t)$ itself, we examine
the spatial integrability of $g_0$.

\begin{remark}\label{rem:iso-tail}
For every datum in Theorem~\ref{thm:main}, $g_0\notin L^2(\R^3)$.
Indeed, $\rho_0\in H^3$ has a $C^1$ representative vanishing at
infinity, so $h_0=\log\rho_0$ tends to $-\infty$ along every ray.
If $Dh_0\in L^2$, polar coordinates and Cauchy--Schwarz would give,
for almost every $\theta\in\mathbb{S}^2$,
\[
 \int_1^\infty |\partial_r h_0(r\theta)|\,dr
 \le \left(\int_1^\infty r^2|Dh_0(r\theta)|^2\,dr\right)^{1/2}
      \left(\int_1^\infty r^{-2}\,dr\right)^{1/2}<\infty .
\]
Then $h_0(r\theta)$ would have a finite limit along each such ray.
It follows from \eqref{eq:reg-shifts} that for $\gamma=1$ and every
$t>0$ the physical velocity $u(t)$ lies outside $L^2$, although
$u(t)+tg_0\to u_0$ in $H^2$ as $t\downarrow0$.
For instance, $\rho_0(x)=\exp(-(1+|x|^2)^{1/8})$ and $u_0=0$
satisfy the hypotheses. Here $g_0=-\tfrac14 x(1+|x|^2)^{-7/8}$ has
magnitude comparable to $|x|^{-3/4}$ at infinity.
The same velocity correction occurs in the field-free equations \cite{GHL}.
\end{remark}

\subsection*{Relation to earlier work}
For the field-free quantum Navier--Stokes equations, J\"ungel used a
logarithmic correction to the velocity to construct periodic weak
solutions and to estimate second derivatives of $\sqrt\rho$
\cite{J10,J11}. The effective velocity also enters the compactness and
periodic existence results of Antonelli and Spirito
\cite{AScompact,ASexist}. Lacroix-Violet and Vasseur obtained global weak
solutions and a semiclassical limit using a velocity-renormalized
formulation \cite{LVV}. In the whole space Antonelli, Hientzsch and
Spirito treated both positive and zero far-field densities for $\gamma>1$
\cite{AHS}. Carles, Carrapatoso and Hillairet constructed whole-space
isothermal weak solutions \cite{CCH}.

Strong solutions of Korteweg models near a stable equilibrium were
constructed by Danchin and Desjardins \cite{DD}. For isothermal quantum
capillarity, Haspot studied continuation through estimates on the
reciprocal density \cite{Haspot}. Huang, Gu and Lei discuss the
dependence on the integrability exponent in the weighted velocity
estimate used there \cite[Remark 1.1]{HGL}. For large data with positive
far-field density, whole-space results include \cite{HGL,HLZ}.
On the torus, Huang, Meng and Zhang proved global regularity for large
positive data \cite{HMZ}. For the quantum coefficients in three
dimensions, Huang and Meng obtained an upper bound on the density
uniform in time for $1\le\gamma<11/4$ \cite[Theorem 1.2]{HM}.
In smooth bounded domains, Li, L\"u and Zhou proved global strong
solvability for $1\le\gamma<3$ in three dimensions under their boundary
curvature condition, with Navier slip for the velocity and a Neumann
condition for the density \cite[Theorem 1.1]{LLZ}.

The closest result for far-field vacuum is due to Gu, Huang and Lei
\cite[Theorems 1.1 and 1.2]{GHL}. For the field-free Cauchy problem
they proved global strong solvability and pointwise comparison with the
initial density. In three dimensions their polytropic theorem assumes
$1<\gamma<8/3$ and $\gamma\rho_0^{\gamma-1}g_0\in L^2$.
Their isothermal theorem permits infinite mass. We use their fixed
density weight and velocity comparisons. In the present problem finite
mass controls the long-range electric field for all pressure laws.
The estimate on the density used to control this field also gives an
integrable comparison coefficient throughout $1<\gamma\le3$.

For Poisson coupling in three dimensions, Yang and Ju proved global
weak existence and an incompressible limit for well-prepared data
\cite{YangJu}. Yang and Li constructed periodic weak solutions with a
singular pressure \cite{YangLi}. For periodic flows with friction and a
doping profile, Antonelli, Cianfarani Carnevale, Lattanzio and Spirito
derived the quantum drift--diffusion limit \cite{ACCLS}. A finite-energy
weak existence result for the high-friction system was obtained in
\cite{Cianfarani}. In one dimension large-data strong solutions are
known for both field-free equations with unequal dispersion and
viscosity parameters \cite{CZ} and a Poisson model away from vacuum
\cite{LiuZhang}. In three dimensions Tong and Xia established global
existence and algebraic decay for small perturbations of a nonconstant
steady state \cite{TongXia}. For Poisson--Korteweg models, Xu and Li
treated a positive background charge under critical Besov smallness
conditions that allow highly oscillatory initial velocities
\cite[Theorem 1.1]{XuLi}, and Gao and Xu obtained global strong
solutions for a class of large data \cite{GaoXu}.

\subsection*{Outline of the proof}
The effective velocity $v=u+D\log\rho$ satisfies a parabolic equation
in which the third-order capillary term cancels. For positive density
the transformed equations read
\begin{equation*}
 (\partial_t-\Delta)\rho=-\Div(\rho v),\qquad
 \Lop v=-\gamma\rho^{\gamma-1}g-E,\qquad
 \Lop h=-|g|^2-\Div v ,
\end{equation*}
where
\begin{equation*}
 \Lop=\partial_t-\Delta+(v-2g)\cdot D .
\end{equation*}
For the smooth approximations of Corollary~\ref{cor:approx-smooth},
the work of the electric field in the effective momentum equation is
\[
 -\int\rho v\cdot E\dx
 =-\frac{d}{dt}\frac{\ell^2}{2}\|E\|_2^2-\ell^{-2}\|\rho\|_2^2 .
\]
The last term follows from $D\rho=\rho g$ and
$\Div E=-\ell^{-2}\rho$ and corresponds to the Poisson contribution in
the periodic entropy inequality of \cite{ACCLS}.
Combining the two momentum balances gives the energy estimate in
Lemma~\ref{lem:energy}, which controls the Fisher information
$\int\rho|g|^2$ and the dissipation $\int_0^S\int\rho|Dv|^2$.
The identity
\[
 \frac{d}{dt}\frac12\int\rho|g|^2\dx=\int\rho Dg:Du\dx
\]
then gives a bound for $\int_0^S\int\rho|D^2\log\rho|^2$
and hence for $\int_0^S\|D^2\sqrt\rho\|_2^2\,dt$.
Applying $\|f\|_\infty^2\lesssim\|Df\|_2\|D^2f\|_2$ to
$f=\sqrt\rho$ yields
\begin{equation}\label{eq:critical-intro}
 \int_0^S\|\rho(t)\|_\infty^2\,dt\le C_T .
\end{equation}
Here $S\le T$ is within the lifespan of a smooth approximation and
$C_T$ is independent of $S$ and of the approximation index, as proved
in Proposition~\ref{prop:critical-density}.
The corresponding density estimate in a bounded domain appears in
\cite[Corollary 3.2]{LLZ}.

Estimate \eqref{eq:critical-intro} has two uses. First, splitting the
Newtonian kernel at a suitable radius gives
$\|E\|_\infty\lesssim\ell^{-2}\|\rho\|_1^{1/3}\|\rho\|_\infty^{2/3}$.
Corollary~\ref{cor:field-time} combines this bound with conservation
of mass to obtain $E\in L^3(0,S;L^\infty)$ before any pointwise
velocity estimate is available.
Second, for $1<\gamma\le3$ the comparison for the smooth
regularization $s+c/\alpha$ of $|v|+(\sqrt\gamma/\alpha)\,\rho^\alpha$,
where $s=(|v|^2+\varepsilon^2)^{1/2}$ and $c=\sqrt\gamma\,\rho^\alpha$
with $\alpha=(\gamma-1)/2$, has
growth coefficient bounded by a multiple of
$\|\rho\|_\infty^{\gamma-1}$. Its time integral is controlled by
\eqref{eq:critical-intro}, with equality of the exponents when $\gamma=3$.

For the isothermal comparison and the continuation argument we use a spatial weight
$\omega=\ee^H$, where $H=\eta*\log\rho_0$ and $\eta$ is a fixed
nonnegative smooth kernel of integral one supported in the unit ball.
The Lipschitz bound on $\log\rho_0$ gives
\begin{equation*}
 \ee^{-L}\rho_0\le\omega\le\ee^L\rho_0,\qquad L=\|g_0\|_\infty .
\end{equation*}
At $\gamma=1$, an upper bound for $\tfrac12|v|^2+h-H$ and a lower
estimate for $h(t)-h(0)$ give a closed integral inequality for
$\|v\|_\infty^2$. For either pressure law, the resulting velocity bound
yields upper and positive lower bounds on $a=\rho/\omega$ and an upper
bound on $|b|=|av|$. Proposition~\ref{prop:relative-continuation} uses these bounds
to continue the local solution for $(a,b)$ past any finite time. Section~\ref{sec:higher} obtains
the Sobolev estimates and strong traces by subtracting smooth profiles
from $v$ and $g$. Section~\ref{sec:proof} recovers the original variables
and proves uniqueness. Appendix~\ref{app:local} provides the local theory
and the propagation of regularity for the approximations.

\subsection*{Notation}
Repeated spatial indices are summed from $1$ to $3$,
$\bx_j$ denotes the $j$th coordinate vector, and $A:B=\sum_{ij}A_{ij}B_{ij}$.
We use Euclidean norms for vectors and Frobenius norms for matrices,
and write $\|\cdot\|_q$ for the spatial $L^q$ norm.
The Fourier transform is unitary with phase $\ee^{-ix\cdot\xi}$, and
$\|f\|_{H^k}=\|\langle\xi\rangle^k\widehat f\|_2$, where
$\langle\xi\rangle=(1+|\xi|^2)^{1/2}$. We write $C_0$ for continuous
functions vanishing at infinity, $\BUC$ for bounded uniformly continuous
functions and $[f]_\theta$ for the spatial H\"older seminorm of order
$\theta$. The projections $\Pi_N$ and $\Pi_{>N}$ restrict Fourier
frequencies to $|\xi|\le N$ and $|\xi|>N$, respectively.
The heat kernel is $G_t(x)=(4\pi t)^{-3/2}\exp(-|x|^2/(4t))$ for $t>0$.
The Newtonian field operator is
\begin{equation*}
 \mathcal N[\varphi]=D\bigl((4\pi\ell^2|\cdot|)^{-1}*\varphi\bigr),
\end{equation*}
where $*$ denotes spatial convolution. We write $\eta_r(x)=r^{-3}\eta(x/r)$
and $\chi_R(x)=\chi(x/R)$, where $\chi\in C_c^\infty(B_2(0))$ satisfies
$0\le\chi\le1$ and $\chi=1$ on $B_1(0)$. In an estimate $A\lesssim B$,
the implicit constant may depend on the dimension, the fixed kernel and
the differentiation order. Other dependencies are indicated when relevant.
The symbol $C_T$ denotes a constant with the dependence specified in
Theorem~\ref{thm:main} and may vary between estimates.

\section{Effective equations and smooth approximations}\label{sec:prelim}

\subsection{The effective system}\label{subsec:effective}

Because $\rho D\log\rho=D\rho$, adding the differentiated continuity
equation to the momentum equation gives an equation for $\rho v$.
In this sum the viscous and capillary stresses combine into a diffusion
term for $v$.

\begin{lemma}\label{lem:effective-system}
Let $(\rho,u,V)$ be a smooth solution of \eqref{eq:model} with $\rho>0$.
Set $h=\log\rho$, $g=Dh$, $v=u+g$, $E=DV$ and
$\Lop=\partial_t-\Delta+(v-2g)\cdot D$. Then
\begin{align}
 (\partial_t-\Delta)\rho&=-\Div(\rho v),\label{eq:heat-rho}\\
 (\partial_t-\Delta)h&=|g|^2-v\cdot g-\Div v,\label{eq:heat-h}\\
 \Lop h&=-|g|^2-\Div v.\label{eq:Lop-h}
\end{align}
The effective velocity and momentum satisfy
\begin{align}
 \Lop v&=-\gamma\rho^{\gamma-1}g-E,\label{eq:effective}\\
 (\partial_t-\Delta)(\rho v_i)
 &=-\partial_j(\rho v_i v_j)-\partial_i\rho^\gamma-\rho E_i .
 \label{eq:heat-m}
\end{align}
Conversely, suppose that smooth $\rho>0$ and $\rho v$ satisfy
\eqref{eq:heat-rho} and \eqref{eq:heat-m} with $E=\mathcal N[\rho]$.
Then $u=v-D\log\rho$ satisfies \eqref{eq:mass-pde} and \eqref{eq:mom}
with this field.
\end{lemma}

\begin{proof}
The continuity equation gives $h_t+u\cdot g=-\Div u$.
The identities $D\rho=\rho g$ and $\Div v=\Div u+\Delta h$ give
\eqref{eq:heat-rho}--\eqref{eq:Lop-h}. Differentiating the continuity
equation in $x_i$ gives
\begin{equation}\label{eq:grad-mass}
 \partial_t(\rho g_i)+\partial_j(\rho u_jg_i)
 =-\partial_j(\rho\,\partial_i u_j).
\end{equation}
Adding this identity to the $i$th momentum equation cancels the
transposed part of the viscous stress. The remaining viscous and
quantum stresses sum to $\partial_j(\rho\,\partial_jv_i)$, since
$\partial_ig_j=\partial_jg_i$. Thus
\begin{equation}\label{eq:material}
 \rho(\partial_t+u\cdot D)v_i
 =\partial_j(\rho\,\partial_jv_i)-\partial_i\rho^\gamma-\rho E_i .
\end{equation}
Division by $\rho$ gives \eqref{eq:effective}.
Using the product rule in $(\partial_t-\Delta)(\rho v_i)$, the terms
$2\rho g\cdot Dv_i$ from \eqref{eq:effective} and $-2D\rho\cdot Dv_i$
cancel. Together with \eqref{eq:heat-rho}, this gives \eqref{eq:heat-m}.

For the converse, substituting $u=v-D\log\rho$ into \eqref{eq:heat-rho}
gives the continuity equation. Expanding \eqref{eq:heat-m} and using
\eqref{eq:heat-rho} together with the positivity of $\rho$ gives
\eqref{eq:material}. Writing the left-hand side of
\eqref{eq:material}, via the continuity equation, as
$\partial_t(\rho v_i)+\partial_j(\rho u_jv_i)$ and subtracting
\eqref{eq:grad-mass} recovers \eqref{eq:mom}.
\end{proof}

\subsection{A fixed weight and the conservative system}\label{subsec:weight}

Let $\eta$ be the mollifier fixed in Section~\ref{sec:introduction}.
The function $h_0=\log\rho_0$ is locally Sobolev by positivity and
has bounded gradient $g_0$. Its continuous representative is
globally Lipschitz with constant $L=\|g_0\|_\infty$. Define
\begin{equation*}
 H=\eta*h_0,\qquad d=DH=\eta*g_0,\qquad \omega=\ee^H,\qquad
 a=\rho/\omega,\qquad b=\rho v/\omega .
\end{equation*}
We use the same weight $\omega$ for every approximation. Its comparison
with $\rho_0$ will translate bounds on $a$ into pointwise bounds on
$\rho$, while its derivatives enter the equations for $a$ and $b$.

\begin{lemma}\label{lem:weight}
The weight satisfies
\begin{equation}\label{eq:weight-bounds}
 |H-h_0|\le L,\qquad \ee^{-L}\rho_0\le\omega\le\ee^L\rho_0,
 \qquad |d|\le L .
\end{equation}
For each $k\ge1$, $\|D^kd\|_\infty\le C_{k,\eta}L$. For every
$q\in\R$ and every integer $k\ge0$,
\begin{equation}\label{eq:weight-powers}
 |D^k(\omega^q)|\le C_{k,q,L,\eta}\,\omega^q .
\end{equation}
\end{lemma}

\begin{proof}
The first inequality follows from
$|H(x)-h_0(x)|\le L\int |y|\eta(y)\,dy\le L$.
Exponentiation gives the bounds for $\omega$, and convolution gives
$d=\eta*g_0$ and $D^kd=D^k\eta*g_0$. Successive differentiation of
$D(\omega^q)=q\omega^qd$ writes each derivative of $\omega^q$ as
$\omega^q$ times a polynomial in $d$ and its derivatives.
\end{proof}

We now divide the equations for $\rho$ and $\rho v$ in
Lemma~\ref{lem:effective-system} by $\omega$. Keeping the fluxes in
divergence form allows us to put the spatial derivative on the heat
kernel in the integral formulation.

\begin{lemma}\label{lem:normalized-system}
Let $\rho=\omega a>0$ and $\rho v=\omega b$ be smooth. The equations
\eqref{eq:heat-rho} and \eqref{eq:heat-m} with $E=\mathcal N[\rho]$
are equivalent to the system for $U=(a,b)$ given by
\begin{equation}\label{eq:U-system}
 U_t-\Delta U=\partial_j\mathcal F_j(U)+\mathcal J(U).
\end{equation}
The fluxes are
\begin{equation}\label{eq:flux}
 \mathcal F_j^a=2d_ja-b_j,\qquad
 \mathcal F_j^{b_i}=2d_jb_i-\frac{b_ib_j}{a}
   -\delta_{ij}\omega^{\gamma-1}a^\gamma
\end{equation}
and the source is
\begin{equation}\label{eq:source}
 \begin{aligned}
 \mathcal J^a&=(|d|^2-\Delta H)a-d\cdot b,\\
 \mathcal J^{b_i}&=(|d|^2-\Delta H)b_i
 -d_j\Bigl(\frac{b_ib_j}{a}+\delta_{ij}\omega^{\gamma-1}a^\gamma\Bigr)
 -a\,\mathcal N[\omega a]_i .
 \end{aligned}
\end{equation}
\end{lemma}

\begin{proof}
Since $d=DH$ and $\omega=\ee^H$, we have
\[
 \omega^{-1}\Delta(\omega y)
 =\Delta y+2d\cdot Dy+(\Div d+|d|^2)y,\qquad
 \omega^{-1}\partial_j(\omega y)=\partial_jy+d_jy .
\]
Apply these identities to \eqref{eq:heat-rho} and \eqref{eq:heat-m}.
The momentum flux divided by $\omega$ is $b_ib_j/a$, and the
corresponding pressure is $\rho^\gamma/\omega=\omega^{\gamma-1}a^\gamma$.
Writing $2d\cdot Dy$ in divergence form gives \eqref{eq:flux} and
\eqref{eq:source}, and hence \eqref{eq:U-system}. Multiplication by
$\omega$ reverses these identities.
\end{proof}

The source in \eqref{eq:source} contains the nonlocal product
$a\,\mathcal N[\omega a]$. For local existence we need a Lipschitz
estimate for this product in the supremum norm. To recover the spatial
derivatives of the field, we also need an estimate in Sobolev spaces.

\begin{lemma}\label{lem:field-reg}
For $\varphi\in L^1\cap L^\infty(\R^3)$,
\begin{equation}\label{eq:weighted-field}
 \|\mathcal N[\varphi]\|_\infty
 \le\ell^{-2}\left(\|\varphi\|_\infty+\frac1{4\pi}\|\varphi\|_1\right).
\end{equation}
If $\varphi\in L^1\cap\BUC$, then $\mathcal N[\varphi]\in\BUC$.
For $a,\tilde a\in\BUC$ with
$\max(\|a\|_\infty,\|\tilde a\|_\infty)\le A$,
\begin{equation}\label{eq:field-lipschitz}
 \|a\,\mathcal N[\omega a]-\tilde a\,\mathcal N[\omega\tilde a]\|_\infty
 \le 2C_\omega A\|a-\tilde a\|_\infty,\qquad
 C_\omega=\ell^{-2}\left(\|\omega\|_\infty+\frac{\|\omega\|_1}{4\pi}\right).
\end{equation}
For every integer $k\ge1$ and $\varphi\in L^1\cap H^{k-1}$,
\begin{equation}\label{eq:field-sobolev}
 \ell^4\|\mathcal N[\varphi]\|_{H^k}^2
 =\int_{\R^3}\frac{\langle\xi\rangle^{2k}}{|\xi|^2}
                |\widehat\varphi(\xi)|^2\,d\xi
 \le C_k\bigl(\|\varphi\|_1^2+\|\varphi\|_{H^{k-1}}^2\bigr).
\end{equation}
\end{lemma}

\begin{proof}
The kernel of $\mathcal N$ is
\begin{equation}\label{eq:newton-kernel}
 K(z)=-\frac{z}{4\pi\ell^2|z|^3}.
\end{equation}
Splitting the convolution at $|z|=1$ and using
$\int_{B_1}|z|^{-2}\,dz=4\pi$ gives \eqref{eq:weighted-field}.
For $\varphi\in L^1\cap\BUC$, translations are continuous in both
norms. Applying \eqref{eq:weighted-field} to differences of translates
gives uniform continuity of the field. In particular
$\|\mathcal N[\omega a]\|_\infty\le C_\omega\|a\|_\infty$, and
$\omega a\in L^1\cap\BUC$. Expanding the difference of the two products
and using linearity gives \eqref{eq:field-lipschitz}.

For \eqref{eq:field-sobolev}, use the Fourier multiplier
$i\xi/(\ell^2|\xi|^2)$. On $|\xi|\le1$, the bound
$|\widehat\varphi|\le(2\pi)^{-3/2}\|\varphi\|_1$ applies and
$|\xi|^{-2}$ is integrable. On $|\xi|>1$, the squared weight
$\langle\xi\rangle^{2k}/|\xi|^2$ is bounded by
$2\langle\xi\rangle^{2(k-1)}$.
\end{proof}

For $a$ bounded away from zero the local terms in
\eqref{eq:flux}--\eqref{eq:source} are locally Lipschitz in $(a,b)$ in
the supremum norm, with constants uniform in $x$. The same holds for
the field term by Lemma~\ref{lem:field-reg}. We use the mild formulation
\begin{equation}\label{eq:mild}
 U(t)=G_t*U_{\mathrm{in}}
 +\int_0^t\Bigl[\partial_jG_{t-s}*\mathcal F_j(U(s))
   +G_{t-s}*\mathcal J(U(s))\Bigr]\,ds ,
\end{equation}
where the free heat flow takes its initial value at $t=0$. Since
$\|DG_t\|_1\lesssim t^{-1/2}$, the time integral of the derivative
kernel is small on short intervals. Together with the Lipschitz bounds
above, this gives a local solution that can be restarted under the
same bounds.

\begin{proposition}\label{prop:local}
Let $U_{\mathrm{in}}=(a_{\mathrm{in}},b_{\mathrm{in}})\in
\BUC(\R^3;\R^4)$ and $\underline a=\inf_xa_{\mathrm{in}}(x)>0$.
There exists $\tau>0$ such that \eqref{eq:mild} has a unique solution
$U\in C([0,\tau];\BUC)$ in the class $\inf_{t,x}a(t,x)>0$. This solution
satisfies $a\ge\underline a/2$. The time $\tau$ may be chosen in terms
of $\|U_{\mathrm{in}}\|_\infty$, $\underline a$, and the fixed
parameters $\ell,\gamma,\eta,L,\|\omega\|_1,\|\omega\|_\infty$.
For two initial data with a common upper norm bound and a common positive
lower bound for the first component, the corresponding solutions with
the same weight satisfy, on a common interval chosen from these bounds,
\begin{equation}\label{eq:stability}
 \sup_{0\le t\le\tau}\|U(t)-\widetilde U(t)\|_\infty
 \le2\|U_{\mathrm{in}}-\widetilde U_{\mathrm{in}}\|_\infty .
\end{equation}
If a solution on $[0,S]$ satisfies
\begin{equation}\label{eq:restart-data}
 \sup_{0\le t\le S}\|U(t)\|_\infty\le A,\qquad
 \inf_{0\le t\le S,\ x\in\R^3}a(t,x)\ge c_0>0 ,
\end{equation}
it extends to $[0,S+\tau_0]$ for some $\tau_0>0$ depending only on
$A,c_0$ and the fixed parameters. If the initial data have bounded
uniformly continuous spatial derivatives through order $m$, then $U$
is continuous on $[0,\tau]$ in the norm given by these derivatives and
the supremum norm.
\end{proposition}

The contraction argument is given in Appendix~\ref{subsec:proof-local}
and spatial regularity in Lemma~\ref{lem:smoothing}. The lower bound in
\eqref{eq:restart-data} applies to $a=\rho/\omega$.

\subsection{Approximation of the initial data}\label{subsec:approx}

We regularize the logarithm of the density and the velocity. For each
integer $n\ge1$, set
\begin{equation*}
 \begin{aligned}
 h_{0,n}&=\eta_{1/n}*h_0,\qquad
 \rho_{0,n}=\ee^{h_{0,n}},\qquad g_{0,n}=\eta_{1/n}*g_0,\\
 u_{0,n}&=\eta_{1/n}*u_0,\qquad
 v_{0,n}=u_{0,n}+g_{0,n},\qquad M_n=\int\rho_{0,n}\,dx .
 \end{aligned}
\end{equation*}
Set $a_{0,n}=\rho_{0,n}/\omega$, $b_{0,n}=a_{0,n}v_{0,n}$, and define
$a_0,b_0$ analogously from $\rho_0,v_0$. Write $M=\int\rho_0\,dx$.

The Lipschitz bound on $h_0$ controls the ratio $\rho_{0,n}/\rho_0$
even in the far field. We also need the normalized data to converge
in $\BUC$ and the differences $g_{0,n}-d$ and $v_{0,n}-\eta*v_0$
to converge in $H^2$.

\begin{lemma}\label{lem:approx-data}
Under the initial hypotheses of Theorem~\ref{thm:main}, the following
properties hold.
\begin{enumerate}
\renewcommand{\theenumi}{\roman{enumi}}
\renewcommand{\labelenumi}{(\theenumi)}
\item\label{itm:density-comparison} $\|h_{0,n}-h_0\|_\infty\le L/n$,
$\ee^{-L/n}\rho_0\le\rho_{0,n}\le\ee^{L/n}\rho_0$, and $M_n\to M$.
\item\label{itm:data-convergence} $g_{0,n}\to g_0$ uniformly, $Dg_{0,n}\to Dg_0$ in $H^1$,
$u_{0,n}\to u_0$ in $H^2$, and $\rho_{0,n}\to\rho_0$ in $L^1\cap H^3$.
\item\label{itm:normalized-data} $\ee^{-2L}\le a_{0,n}\le\ee^{2L}$,
$(a_{0,n},b_{0,n})\to(a_0,b_0)$ in $\BUC$, and the spatial
$1/2$-H\"older seminorms of the pairs are uniformly bounded.
Moreover, $Da_{0,n}\to Da_0$ uniformly and
$\sup_n\|Da_{0,n}\|_{H^2}<\infty$.
\item\label{itm:initial-pressure} If $\gamma>1$, then
$p_{0,n}=\gamma\rho_{0,n}^{\gamma-1}g_{0,n}\to
p_0=\gamma\rho_0^{\gamma-1}g_0$ in $L^2$; for $1<\gamma<5/4$ this uses
the additional hypothesis \eqref{eq:initial-pressure} of
Theorem~\ref{thm:main}, while for $5/4\le\gamma\le3$ it follows from
Lemma~\ref{lem:initial-buc}.
\item\label{itm:shifted-data} The initial values
\begin{equation*}
 w_n(0)=v_{0,n}-\eta*v_0,\qquad z_n(0)=g_{0,n}-d
\end{equation*}
converge in $H^2$ to $v_0-\eta*v_0$ and $g_0-d$, respectively.
\end{enumerate}
\end{lemma}

\begin{proof}
The Lipschitz bound for $h_0$ gives (\ref{itm:density-comparison}) and, for $1\le q\le\infty$,
$\|\rho_{0,n}-\rho_0\|_q\le(\ee^{L/n}-1)\|\rho_0\|_q$.
Lemma~\ref{lem:initial-buc} gives $g_0\in\BUC$. Mollification therefore
yields uniform convergence of $g_{0,n}$ and strong convergence of
$Dg_{0,n},u_{0,n}$ in the indicated Sobolev norms. To obtain
$\rho_{0,n}\to\rho_0$ in $H^3$, differentiate $D\rho_{0,n}=\rho_{0,n}g_{0,n}$
and the corresponding identity for $\rho_0$. Up to third order the
resulting terms are products of $\rho$ with powers of $g$, with $Dg$,
with $g$ and $Dg$, or with $D^2g$. In products containing $Dg$ or
$D^2g$, the differentiated factor converges in $L^2$ and the remaining
factors converge uniformly. In the other products use the $L^2$
convergence of the density and uniform convergence of $g_{0,n}$. This
proves (\ref{itm:data-convergence}).

For $\psi\in W^{1,2}_{\mathrm{loc}}$ with $D\psi\in L^2$, the
translation estimate and Minkowski's inequality give
\begin{equation}\label{eq:translation}
 \|\eta_r*\psi-\psi\|_2
 \le\int\eta_r(y)|y|\|D\psi\|_2\,dy\le r\|D\psi\|_2 .
\end{equation}
Indeed, $\psi(x-y)-\psi(x)$ is the integral of
$-y\cdot D\psi(x-sy)$ over $0\le s\le1$ for almost every $(x,y)$.
In particular $q_n=g_{0,n}-d$ is uniformly bounded in $L^2\cap L^\infty$
and $Dq_n$ is uniformly bounded in $H^1$.

Part (\ref{itm:density-comparison}) and Lemma~\ref{lem:weight} give the bounds on $a_{0,n}$.
Since $a_{0,n}=a_0\exp(h_{0,n}-h_0)$ and $a_0$ is bounded,
part (\ref{itm:density-comparison}) also gives $a_{0,n}\to a_0$ uniformly.
Together with (\ref{itm:data-convergence}), this gives uniform convergence of $b_{0,n}=a_{0,n}v_{0,n}$
and $Da_{0,n}=a_{0,n}q_n$. The uniform H\"older bounds follow from
Lemma~\ref{lem:initial-buc}. Differentiating $Da_{0,n}=a_{0,n}q_n$
twice produces terms bounded in modulus by
$C a_{0,n}(|q_n|^3+|q_n||Dq_n|+|D^2q_n|)$.
The corresponding first derivative is bounded by
$C a_{0,n}(|q_n|^2+|Dq_n|)$. The bounds on $q_n$ give (\ref{itm:normalized-data}).

For $\beta=\gamma-1>0$, the assumptions of the theorem and
Lemma~\ref{lem:initial-buc} give $p_0\in L^2$. Splitting the difference
of the products and applying (\ref{itm:density-comparison}) and \eqref{eq:translation}, we obtain
\[
 \|p_{0,n}-p_0\|_2\le
 \gamma\|\rho_{0,n}\|_\infty^\beta n^{-1}\|Dg_0\|_2
 +(\ee^{\beta L/n}-1)\|p_0\|_2 .
\]
This proves (\ref{itm:initial-pressure}). Finally, \eqref{eq:translation} shows both that
$g_0-d\in L^2$ and that $\|g_{0,n}-g_0\|_2\to0$. Combining this with the
convergence of $Dg_{0,n}$ in $H^1$ proves (\ref{itm:shifted-data}) for $z_n(0)$.
The assertion for $w_n(0)$ follows by adding
$u_{0,n}-\eta*u_0\to u_0-\eta*u_0$ in $H^2$.
\end{proof}

The lower bound on $a_{0,n}$ allows us to apply
Proposition~\ref{prop:local} to each pair of mollified data. The
resulting solutions provide the smooth fields on which we derive the
a priori estimates.

\begin{corollary}\label{cor:approx-smooth}
Let $[0,S_n)$ be the maximal interval on which the solution $U_n$ with
data $(a_{0,n},b_{0,n})$ exists in the positive class of
Proposition~\ref{prop:local}. For every $S<S_n$ and every integer $k\ge0$,
\begin{align*}
 &\sup_{0\le t\le S}\|D^kU_n(t)\|_\infty<\infty,\\
 &\inf_{[0,S]\times\R^3}a_n>0.
\end{align*}
The fields $\rho_n=\omega a_n$, $v_n=b_n/a_n$,
$g_n=d+Da_n/a_n$ and $u_n=v_n-g_n$ are smooth and satisfy the original
system with $E_n=\mathcal N[\rho_n]$. Their spatial derivatives of every
order are bounded on $[0,S]\times\R^3$.
\end{corollary}

\begin{proof}
Mollification gives bounded derivatives of every order for $g_{0,n}$
and $u_{0,n}$, and repeated differentiation of $Da_{0,n}=a_{0,n}(g_{0,n}-d)$
gives the same property for $a_{0,n},b_{0,n}$. Boundedness at the next
order gives uniform continuity at each order. Proposition~\ref{prop:local}
and Lemma~\ref{lem:smoothing} apply successively on every compact
subinterval of the maximal lifespan. There $a_n$ has a positive lower
bound, and Lemma~\ref{lem:weight} bounds each derivative of $\omega$.
Thus the stated fields and their spatial derivatives are bounded.
The equation gives smoothness in time, with one-sided derivatives
at $t=0$. Lemma~\ref{lem:normalized-system} gives
\eqref{eq:heat-rho} and \eqref{eq:heat-m}, and
Lemma~\ref{lem:effective-system} recovers the original equations.
\end{proof}

\section{Energy identities and the critical density integral}\label{sec:energy}

Fix $T>0$, $n\ge1$ and $0<S<S_n$ with $S\le T$. Throughout this section
we suppress the approximation index except in the mass $M_n$. On
$[0,S]$, Corollary~\ref{cor:approx-smooth} gives bounds for $u,v,g$ and
their spatial derivatives, and $\rho\le C_{n,S}\omega$. In particular
$\rho$ and $\rho^\gamma$ are dominated by integrable functions.
We use these bounds to justify identities for a fixed approximation.
The estimates denoted by $C_T$ below are independent of $n$ and $S$.

\subsection{Energy identities}\label{subsec:energy-identities}

Testing the two momentum equations with $u$ and $v$ produces different
electric work terms. To combine the resulting balances, we also track
the capillary work through the Fisher information
$I(t)=\tfrac12\int\rho|g|^2\,dx$.

\begin{lemma}\label{lem:energy-identities}
For every $0\le t\le S$, $\int\rho(t)\,dx=M_n$. For almost every
$t\in(0,S)$ the physical and effective electric work satisfy
\begin{align}
 \frac d{dt}\frac{\ell^2}2\|E\|_2^2
 &=\int\rho u\cdot E\,dx ,\label{eq:field-work}\\
 \int\rho v\cdot E\,dx
 &=\frac d{dt}\frac{\ell^2}2\|E\|_2^2+\ell^{-2}\|\rho\|_2^2,
 \label{eq:effective-work}
\end{align}
and the Fisher information obeys
\begin{equation}\label{eq:fisher-work}
 I'(t)=\int\rho Dg:Du\,dx .
\end{equation}
For almost every $t\in(0,S)$ the kinetic energies satisfy
\begin{align*}
 \frac d{dt}\frac12\int\rho|u|^2\,dx
 +2\int\rho|\dsym u|^2\,dx
 &=\int\rho^\gamma\Div u\,dx-\int\rho Dg:Du\,dx
       -\int\rho u\cdot E\,dx ,\\
 \frac d{dt}\frac12\int\rho|v|^2\,dx+\int\rho|Dv|^2\,dx
 &=\int\rho^\gamma\Div v\,dx-\int\rho v\cdot E\,dx .
\end{align*}
\end{lemma}

\begin{proof}
Use the cutoff $\chi_R$ of Section~\ref{sec:introduction}. The boundary
terms arising from the continuity and momentum equations have the form
\begin{equation}\label{eq:cutoff-flux}
 \left|\int J\cdot D\chi_R\,dx\right|
 \le R^{-1}\|D\chi\|_\infty\|J\|_1 .
\end{equation}
The corresponding fluxes are products of $\rho$ or $\rho^\gamma$ with
bounded factors. In particular the continuity equation tested against
$\chi_R$ yields $\int\rho(t)\,dx=M_n$ by time integration and
dominated convergence. Lemma~\ref{lem:approx-data} also gives
$\ee^{-L/n}M\le M_n\le\ee^{L/n}M$.

The Fourier representation of the field is
$\widehat E=i\xi\widehat\rho/(\ell^2|\xi|^2)$ for $\xi\ne0$.
Splitting at $|\xi|=1$ and using
$|\widehat\rho|\le(2\pi)^{-3/2}M_n$ gives
\begin{equation}\label{eq:field-low}
 \ell^4\|E\|_2^2
 =\int\frac{|\widehat\rho(\xi)|^2}{|\xi|^2}\,d\xi
 \le\frac{4\pi M_n^2}{(2\pi)^3}+\|\rho\|_2^2 .
\end{equation}
The same estimate for differences, with their $L^1$ norms in place
of $M_n$, gives $E\in C([0,S];L^2)$ from
$\rho=\omega a\in C([0,S];L^1\cap L^2)$. Let $P_\parallel$ denote the
orthogonal projection in $L^2(\R^3;\R^3)$ with Fourier symbol
$\xi\otimes\xi/|\xi|^2$. The continuity equation implies
$\ell^2 E_t=P_\parallel(\rho u)$ in distributions.
As $\rho u\in L^\infty(0,S;L^2)$, the field belongs to
$W^{1,\infty}(0,S;L^2)$. The Hilbert-space chain rule,
self-adjointness of $P_\parallel$ and $P_\parallel E=E$ yield
\eqref{eq:field-work}. To obtain \eqref{eq:effective-work}, substitute $v=u+g$ and integrate
$D\rho=\rho g$ by parts, using $\Div E=-\ell^{-2}\rho$.
The flux $\rho E$ belongs to $L^1$ by
\eqref{eq:weighted-field}.

For a smooth scalar $\psi$, the continuity equation gives the identity
\[
 \psi\bigl[(\rho\psi)_t+\Div(\rho u\psi)\bigr]
 =\tfrac12\partial_t(\rho\psi^2)+\tfrac12\Div(\rho u\psi^2).
\]
Testing \eqref{eq:grad-mass} by $\chi_R g_i$, summing over $i$ and
integrating by parts yields
\[
 \frac d{dt}\frac12\int\chi_R\rho|g|^2\,dx
 =\int\chi_R\rho\,\partial_i u_j\,\partial_jg_i\,dx
 +\int\left(\tfrac12\rho u_j|g|^2+\rho g_i\partial_iu_j\right)
          \partial_j\chi_R\,dx .
\]
Integrating in time and using \eqref{eq:cutoff-flux} and the symmetry
of $Dg$ proves \eqref{eq:fisher-work}.
Testing the physical momentum equation by $\chi_R u$ and
\eqref{eq:material} by $\chi_R v$ gives the two kinetic-energy
balances as $R\to\infty$. The fluxes on the cutoff annulus are
bounded in modulus by $C_{n,S}(\rho+\rho^\gamma)$, so the identities
hold in integrated form on every subinterval of $[0,S]$.
\end{proof}

The electric and Fisher identities in Lemma~\ref{lem:energy-identities}
account for the field and capillary work in the sum of the kinetic
balances. For $\gamma>1$ the pressure contributes an internal energy
term, while for $\gamma=1$ we estimate its work directly. With
$\sigma=\mathbf 1_{\{\gamma=1\}}$, define
\begin{equation*}
 \mathscr E(t)=\frac12\int\rho\bigl(|u|^2+|g|^2+|v|^2\bigr)\,dx
 +\ell^2\|E\|_2^2
 +\begin{cases}\displaystyle\frac{2}{\gamma-1}\int\rho^\gamma\,dx,
                       &\gamma>1,\\[1mm] 0,&\gamma=1 .\end{cases}
\end{equation*}

\begin{lemma}\label{lem:energy}
For $0\le t\le S$,
\begin{multline}\label{eq:energy-bound}
 \mathscr E(t)
 +\int_0^t\int\rho\bigl((2-\sigma)|\dsym u|^2+|Dv|^2\bigr)\,dx\,ds\\
 +\int_0^t\left(\frac4\gamma\|D\rho^{\gamma/2}(s)\|_2^2
              +\ell^{-2}\|\rho(s)\|_2^2\right)\,ds
 \le\mathscr E(0)+3\sigma M_nt .
\end{multline}
In particular,
\begin{equation}\label{eq:energy-uniform}
 \sup_{0\le t\le S}\mathscr E(t)
 +\int_0^S\!\int\rho\bigl((2-\sigma)|\dsym u|^2+|Dv|^2\bigr)\,dx\,ds
 +\int_0^S\|D\rho^{\gamma/2}(s)\|_2^2\,ds\le C_T ,
\end{equation}
with $C_T$ independent of $n$ and $S$.
\end{lemma}

\begin{proof}
Using $v=u+g$ and integrating by parts gives
\begin{equation}\label{eq:pressure-work}
 \int\rho^\gamma\Div v\,dx=\int\rho^\gamma\Div u\,dx
 -\gamma\int\rho^\gamma|g|^2\,dx
 =\int\rho^\gamma\Div u\,dx-\frac4\gamma\|D\rho^{\gamma/2}\|_2^2 .
\end{equation}
For $\gamma>1$ the continuity equation gives
$\partial_t\rho^\gamma+\Div(\rho^\gamma u)
=-(\gamma-1)\rho^\gamma\Div u$. For $\gamma=1$ we estimate the
pressure work directly, without imposing an integrability condition on
$\rho\log\rho$. Adding the two kinetic-energy balances from Lemma~\ref{lem:energy-identities} and
using \eqref{eq:field-work}, \eqref{eq:effective-work},
\eqref{eq:fisher-work} and \eqref{eq:pressure-work} gives
\begin{multline*}
 \frac d{dt}\mathscr E
 +2\int\rho|\dsym u|^2\,dx+\int\rho|Dv|^2\,dx
 +\frac4\gamma\|D\rho^{\gamma/2}\|_2^2+\ell^{-2}\|\rho\|_2^2
 =2\sigma\int\rho\Div u\,dx .
\end{multline*}
In the isothermal case,
$|\Div u|\le\sqrt3|\dsym u|$ and Cauchy--Schwarz give
\[
 \left|2\int\rho\Div u\,dx\right|
 \le2\sqrt{3M_n}\left(\int\rho|\dsym u|^2\,dx\right)^{1/2}
 \le3M_n+\int\rho|\dsym u|^2\,dx .
\]
Integrating in time gives \eqref{eq:energy-bound}.

The initial kinetic and Fisher energies are bounded uniformly by
$C M_n(\|u_0\|_\infty^2+L^2)$. For $\gamma>1$ the initial pressure
integral is at most $M_n\|\rho_{0,n}\|_\infty^{\gamma-1}$, and
\eqref{eq:field-low} bounds the initial field energy. Thus
Lemma~\ref{lem:approx-data} gives $\sup_n\mathscr E(0)<\infty$, with a
bound depending only on the data and fixed parameters. Since $M_n\le\ee^{L/n}M$ and $S\le T$,
\eqref{eq:energy-bound} implies \eqref{eq:energy-uniform}.
\end{proof}

\subsection{The critical density integral}\label{subsec:critical}

The dissipation in Lemma~\ref{lem:energy} controls $Dv$ with weight
$\rho$. In the expansion of $|Dv|^2=|Du+Dg|^2$, the Fisher identity
\eqref{eq:fisher-work} accounts for the cross term. The resulting
bound on $\int_0^S\!\int\rho|Dg|^2\,dx\,dt$ controls $D^2\sqrt{\rho}$ in
$L^2(0,S;L^2)$, which is the derivative estimate needed below.

\begin{proposition}\label{prop:critical-density}
For $f=\sqrt{\rho}$, the smooth approximation on $[0,S]$ satisfies
\begin{equation}\label{eq:Dg-bound}
 \int_0^S\int\rho |Dg|^2\,dx\,ds\le C_T
\end{equation}
and
\begin{equation}\label{eq:critical}
 \int_0^S\|\rho(s)\|_\infty^2\,ds
 \le C\sup_{0\le s\le S}\|Df(s)\|_2^2\int_0^S\|D^2f(s)\|_2^2\,ds
 \le C_T ,
\end{equation}
where $C_T$ is independent of $n$ and $S$.
\end{proposition}

\begin{proof}
Expanding $Dv=Du+Dg$ and integrating the cross term by
\eqref{eq:fisher-work}, we obtain
\begin{equation*}
 \int_0^t\int\rho|Dv|^2\,dx\,ds
 =\int_0^t\int\rho\bigl(|Du|^2+|Dg|^2\bigr)\,dx\,ds+2I(t)-2I(0).
\end{equation*}
Since $2I(0)\le L^2 M_n$, Lemma~\ref{lem:energy} gives
\eqref{eq:Dg-bound}.
At each fixed time we use $D\rho=\rho g$ and integrate by parts to obtain
\[
 \int\rho|g|^4\,dx=-\int\rho\bigl(2g_ig_j\partial_jg_i+
                   |g|^2\Div g\bigr)\,dx
 \le(2+\sqrt3)
    \left(\int\rho|Dg|^2\,dx\right)^{1/2}
    \left(\int\rho|g|^4\,dx\right)^{1/2}.
\]
The boundary flux $\rho g|g|^2$ satisfies \eqref{eq:cutoff-flux}, and
all the integrals are finite for the fixed approximation. Consequently
\begin{equation*}
 \int\rho|g|^4\,dx\le(2+\sqrt3)^2\int\rho|Dg|^2\,dx .
\end{equation*}
For $f=\sqrt\rho$ we have $Df=fg/2$ and
$D^2f=fDg/2+f(g\otimes g)/4$. Hence
\begin{equation}\label{eq:sqrt-hessian}
 \|Df(t)\|_2^2=I(t)/2,\qquad
 \|D^2 f(t)\|_2\le\frac{4+\sqrt3}{4}
           \left(\int\rho|Dg|^2\,dx\right)^{1/2}.
\end{equation}
In particular $f(t)\in H^2$, since $\|f(t)\|_2^2=M_n$.

The three-dimensional interpolation estimate needed here is
\begin{equation}\label{eq:fourier-opt}
 \|f\|_\infty^2\le C\|Df\|_2\|D^2f\|_2 .
\end{equation}
Indeed, splitting $\|\widehat f\|_1$ at $|\xi|=r$ and applying
Cauchy--Schwarz with weights $|\xi|$ and $|\xi|^2$, respectively, bounds
$\|f\|_\infty$ by
$C(r^{1/2}\|Df\|_2+r^{-1/2}\|D^2f\|_2)$. Minimizing over $r>0$ gives
\eqref{eq:fourier-opt}. Squaring it with $f=\sqrt\rho$ and integrating in time gives \eqref{eq:critical}
by \eqref{eq:energy-uniform}, \eqref{eq:Dg-bound} and
\eqref{eq:sqrt-hessian}.
\end{proof}

In the Newtonian field integral, the density supremum controls the
contribution near the singularity and the conserved mass controls the
rest. Balancing the two contributions converts \eqref{eq:critical}
into a time-integrability estimate for the force.

\begin{corollary}\label{cor:field-time}
For the same smooth approximation and every $1\le q\le3$,
\begin{equation}\label{eq:field-time}
 \int_0^S\|E(s)\|_\infty^q\,ds
 \le T^{1-q/3}\left(\int_0^S\|E(s)\|_\infty^3\,ds\right)^{q/3}
 \le C_T ,
\end{equation}
where $C_T$ is independent of $n$ and $S$.
\end{corollary}

\begin{proof}
Split the kernel \eqref{eq:newton-kernel} at radius $r$. For every $r>0$,
\[
 |E(t,x)|\le\ell^{-2}
       \left(r\|\rho(t)\|_\infty+\frac{M_n}{4\pi r^2}\right).
\]
Taking $r=(M_n/(4\pi\|\rho(t)\|_\infty))^{1/3}$ gives
\begin{equation*}
 \|E(t)\|_\infty\le
 \frac{2M_n^{1/3}}{(4\pi)^{1/3}\ell^2}\|\rho(t)\|_\infty^{2/3} .
\end{equation*}
Cubing this bound and integrating, we obtain from
Proposition~\ref{prop:critical-density}
\begin{equation*}
 \int_0^S\|E(s)\|_\infty^3\,ds
 \le\frac{2M_n}{\pi\ell^6}\int_0^S\|\rho(s)\|_\infty^2\,ds\le C_T ,
\end{equation*}
where we also used $M_n\le\ee^{L/n}M$. H\"older's inequality in time
and $S\le T$ now give \eqref{eq:field-time}.
\end{proof}

In particular the source bounds for the velocity comparisons in
Section~\ref{sec:velocity} follow from mass and energy alone.

\section{Bounds for the effective velocity and the relative density}\label{sec:velocity}

Fix $T>0$. Until the continuation argument below, all estimates are on
$[0,S]$ with $S\le T$ and $S<S_n$, where $[0,S_n)$ is the maximal
interval of the $n$-th approximation from Corollary~\ref{cor:approx-smooth}.
We suppress $n$ when working with
one approximation. By \eqref{eq:critical} and \eqref{eq:field-time},
\begin{equation}\label{eq:CT-energy}
 \int_0^S\|\rho(s)\|_\infty^2\,ds
 +\sum_{q=1}^3\int_0^S\|E(s)\|_\infty^q\,ds\le C_T .
\end{equation}
All estimates denoted by $C_T$ are independent of $n$ and $S$.

To use \eqref{eq:CT-energy} in a pointwise velocity bound, we first
need a comparison principle on $\R^3$. The scalar quantities to which
we apply it are bounded above, and their drift is bounded on each
fixed approximation interval.

\begin{lemma}\label{lem:comparison}
Suppose that $Y\in C([0,S]\times\R^3)\cap C^{1,2}((0,S)\times\R^3)$
is bounded above, $\Upsilon$ is a bounded drift, and
$(\partial_t-\Delta+\Upsilon\cdot D)Y\le f(t)$ with $f\in C([0,S])$.
Then
\[
 \sup_xY(t,x)\le\sup_xY(0,x)+\int_0^t f(s)\,ds,\qquad 0\le t\le S .
\]
\end{lemma}

\begin{proof}
Set $\kappa=7+\|\Upsilon\|_\infty^2$ and
$\Theta=\ee^{\kappa t}(1+|x|^2)$. Since
$2\Upsilon\cdot x\ge-\|\Upsilon\|_\infty^2-|x|^2$,
$(\partial_t-\Delta+\Upsilon\cdot D)\Theta\ge\ee^{\kappa t}$.
For $\delta>0$, subtract $\sup_xY(0,x)+\int_0^t f+\delta\Theta$ from
$Y$. The difference is initially negative and tends to $-\infty$
uniformly at spatial infinity. At a positive maximum on $[0,t]\times\R^3$,
$t<S$, its parabolic operator would be nonnegative, contrary to the
strict differential inequality. Let $\delta\downarrow0$, then let
$t\uparrow S$ at each fixed spatial point.
\end{proof}

For the logarithmic density, $h=H+\log a$ has at most linear growth
on each fixed approximation interval, since $H$ is Lipschitz and
$\log a$ is bounded there. A heat-kernel representation under this
growth condition will let us compare $h(t)$ with $h(0)$.

\begin{lemma}\label{lem:heat-duhamel}
Let $y$ be continuous on $[0,S]\times\R^3$ with
$|y(t,x)|\le C(1+|x|)$. Suppose that
$(\partial_t-\Delta)y=\partial_jF_j+J$ in distributions on
$(0,S)\times\R^3$, with bounded continuous $F_j,J$. Then, for $t>0$,
\begin{equation}\label{eq:heat-representation}
 y(t)=G_t*y(0)+\int_0^t
 \bigl[\partial_jG_{t-s}*F_j(s)+G_{t-s}*J(s)\bigr]\,ds .
\end{equation}
\end{lemma}

\begin{proof}
For $\psi\in C_c^\infty(\R^3)$ test on $[0,t]$ with
$\chi_R(G_{t-s}*\psi)$, using smooth time cutoffs to recover the endpoint
pairings. The backward heat equation cancels the interior terms involving
$y$. The family $G_r*\psi$, $0\le r\le S$, is bounded in the Schwartz
space, so the linear growth of $y$ and the bounds on $F_j,J$ allow the
cutoff errors to tend to zero as $R\to\infty$. The resulting identity is
\eqref{eq:heat-representation} tested against $\psi$.
The time singularity $\|DG_{t-s}\|_1\lesssim(t-s)^{-1/2}$ is integrable,
and both sides are continuous in space.
\end{proof}

\subsection{The polytropic case}\label{subsec:poly}

For $\gamma>1$, the pressure term in the equation for $v$ suggests
estimating $|v|$ together with $\rho^{(\gamma-1)/2}$. The resulting
comparison involves $\|\rho\|_\infty^{\gamma-1}$, whose time integral
is controlled by \eqref{eq:CT-energy} up to $\gamma=3$.

\begin{proposition}\label{prop:poly-velocity}
For $1<\gamma\le3$, every smooth approximation satisfies
\begin{equation}\label{eq:poly-conclusion}
 \sup_{0\le t\le S}
 \left(\|v(t)\|_\infty+\|\rho(t)\|_\infty^{(\gamma-1)/2}\right)\le C_T .
\end{equation}
\end{proposition}

\begin{proof}
Write $\alpha=(\gamma-1)/2$ and
$c=\sqrt\gamma\,\rho^\alpha$. For $0<\varepsilon\le1$ set
$s=(|v|^2+\varepsilon^2)^{1/2}$, $e=v/s$ and $P=I-e\otimes e$.
We adapt the scalar comparison in
\cite[Section 4]{GHL} to the electric force. The chain rule,
\eqref{eq:Lop-h} and $Dc=\alpha cg$ give
\[
 \Lop c=-\alpha c\Div v-\alpha(1+\alpha)c|g|^2 .
\]
Since $\partial_js=e\cdot\partial_jv$ and $\partial_je=P\partial_jv/s$,
the effective momentum equation gives
\[
 \Lop s=-c^2e\cdot g-e\cdot E
         -s^{-1}\sum_j(\partial_jv)^TP\partial_jv .
\]
In the equation for $s+c/\alpha$, the additional drift $ce$ contributes
$ce\cdot D(c/\alpha)=c^2e\cdot g$, while
$ce\cdot Ds-c\Div v=-c\sum_j\bx_j^TP\partial_jv$.
Completing the square gives
\begin{align*}
 (\Lop+ce\cdot D)(s+c/\alpha)
 & =-(1+\alpha)c|g|^2-e\cdot E
      +\tfrac14 sc^2\mtrace P \\
 &\quad-\frac1s\sum_j
   \left|P^{1/2}\left(\partial_jv+\tfrac12sc\,\bx_j\right)\right|^2 .
\end{align*}
Here $P$ is positive semidefinite and
$\mtrace P=2+\varepsilon^2/s^2\le3$. Thus
\begin{equation*}
 (\Lop+ce\cdot D)(s+c/\alpha)
 \le \|E(t)\|_\infty+
   \tfrac{3\gamma}{4}\|\rho(t)\|_\infty^{\gamma-1}(s+c/\alpha).
\end{equation*}
Set $\Psi(t)=\tfrac{3\gamma}{4}\int_0^t\|\rho(r)\|_\infty^{\gamma-1}\,dr$.
For each fixed approximation, Corollary~\ref{cor:approx-smooth} bounds
the drift $v-2g+ce$ and the function $s+c/\alpha$ on $[0,S]\times\R^3$.
The maps $t\mapsto\rho(t),E(t)$ are continuous into $\BUC$ by
$\rho=\omega a$ and \eqref{eq:weighted-field}.
Applying Lemma~\ref{lem:comparison} to $\ee^{-\Psi}(s+c/\alpha)$ yields
\[
 \|(s+c/\alpha)(t)\|_\infty
 \le \ee^{\Psi(t)}
 \left(\|(s+c/\alpha)(0)\|_\infty+\int_0^t\|E(r)\|_\infty\,dr\right).
\]
The initial norm is bounded independently of $n$ and
$0<\varepsilon\le1$ by Lemma~\ref{lem:approx-data}.
For $1<\gamma<3$, H\"older's inequality gives
\[
 \int_0^S\|\rho\|_\infty^{\gamma-1}\,dt
 \le S^{(3-\gamma)/2}
       \left(\int_0^S\|\rho\|_\infty^2\,dt\right)^{(\gamma-1)/2}.
\]
At $\gamma=3$ the integral on the left is exactly the density integral
in \eqref{eq:CT-energy}. Thus $\sup_{[0,S]}\Psi\le C_T$ throughout the
stated range. Letting $\varepsilon\downarrow0$ gives \eqref{eq:poly-conclusion}.
\end{proof}

\subsection{The isothermal case}\label{subsec:iso}

At $\gamma=1$, the power of $\rho$ used in the polytropic comparison
reduces to a constant. We instead use the change in $h$ from its initial
value. Its lower bound comes from the mass equation and therefore
applies to both pressure laws.

\begin{lemma}\label{lem:log-lower}
For $1\le\gamma\le3$, set $W(t)=\sup_{0\le s\le t}\|v(s)\|_\infty$ and
$L=\|g_0\|_\infty$. Each smooth approximation satisfies, for $0\le t\le S$
and $x\in\R^3$,
\begin{equation}\label{eq:log-lower}
 h(t,x)-h(0,x)
 \ge -4L\sqrt{t/\pi}-\tfrac14\int_0^t\|v(s)\|_\infty^2\,ds
       -4\sqrt{t/\pi}\,W(t).
\end{equation}
\end{lemma}

\begin{proof}
Since $h=H+\log a$, the fixed-approximation bounds on $a,1/a$
and the Lipschitz bound on $H$ give $|h(t,x)|\le C_{n,S}(1+|x|)$.
Apply Lemma~\ref{lem:heat-duhamel} to \eqref{eq:heat-h}, with
$F=-v$ and $J=|g-v/2|^2-|v|^2/4$. The Gaussian moments are
\begin{equation*}
 \int |y|G_r(y)\,dy=4\sqrt{r/\pi},\qquad
 \|DG_r\|_1=2/\sqrt{\pi r},\qquad
 \int |y|^2G_r(y)\,dy=6r .
\end{equation*}
Since $h(0)=h_{0,n}$ is $L$-Lipschitz, the first moment bounds
$G_t*h(0)-h(0)$ from below by $-4L\sqrt{t/\pi}$. The term with $F$
is bounded in modulus by $W(t)\int_0^t\|DG_{t-s}\|_1\,ds$.
Using $J\ge-|v|^2/4$ in \eqref{eq:heat-representation} gives
\eqref{eq:log-lower}.
\end{proof}

For isothermal pressure, $\tfrac12|v|^2+h-H$ satisfies a scalar
parabolic inequality. To extract a bound on $|v|$ from it, we use
Lemma~\ref{lem:log-lower} to control how far $h$ can fall below its
initial value.

\begin{proposition}\label{prop:iso-velocity}
For $\gamma=1$, every smooth approximation satisfies
\begin{equation}\label{eq:iso-conclusion}
 \sup_{0\le t\le S}\|v(t)\|_\infty\le C_T .
\end{equation}
\end{proposition}

\begin{proof}
Following \cite[Section 5]{GHL}, set
$Y=\tfrac12|v|^2+h-H$. As $\Lop v=-g-E$,
$\Lop h=-|g|^2-\Div v$ and $DH=d$, direct calculation gives
\begin{align*}
 \Lop Y
 &=-|Dv+I/2|^2-|g+v/2-d|^2+\tfrac34+\tfrac14|v|^2 \\
 &\quad -2v\cdot d+|d|^2+\Delta H-v\cdot E .
\end{align*}
Young's inequality bounds the two mixed terms by
$\tfrac14|v|^2+8|d|^2+2|E|^2$. Since $h-H=\log a$, $Y$ is bounded
above on the fixed interval. Applying Lemma~\ref{lem:comparison} with
drift $v-2g$ yields
\[
 \sup_xY(t,x)\le\sup_xY(0,x)
   +\bigl(\tfrac34+9L^2+\|\Delta H\|_\infty\bigr)t
   +\tfrac12\int_0^t\|v\|_\infty^2\,ds
   +2\int_0^t\|E\|_\infty^2\,ds .
\]
The bounds $|h(0)-H|\le2L$ and
$\|v(0)\|_\infty\le\|v_0\|_\infty$ control the initial terms.
Write $\tfrac12|v|^2=Y+(H-h(0))+(h(0)-h)$ and take the supremum
up to time $t$. Lemma~\ref{lem:log-lower} and \eqref{eq:CT-energy}
then give, for $W(t)=\sup_{0\le s\le t}\|v(s)\|_\infty$,
\[
 \tfrac12W(t)^2\le C_T+\tfrac34\int_0^t W(s)^2\,ds
                   +4\sqrt{t/\pi}\,W(t).
\]
The last term is at most $\tfrac14W(t)^2+16t/\pi$.
Gr\"onwall's inequality therefore gives \eqref{eq:iso-conclusion}.
\end{proof}

\subsection{Relative density, logarithmic gradient and convergence}\label{subsec:relative}

Both pressure laws now give a bound on $v$ up to any time $S<S_n$
with $S\le T$. To apply the continuation criterion
\eqref{eq:restart-data}, it remains to control $a=\rho/\omega$ from
above and away from zero, since $b=av$.

\begin{proposition}\label{prop:relative-continuation}
For $1\le\gamma\le3$ and every $T>0$, each smooth approximation extends
through $[0,T]$. There exist $c_T>0$ and $A_T<\infty$, independent of $n$,
such that
\begin{equation}\label{eq:U-bounds}
 a\ge c_T,\qquad \|U\|_{L^\infty([0,T]\times\R^3)}\le A_T .
\end{equation}
Moreover,
\[
 \sup_{0\le t\le T}
 \bigl([U(t)]_{1/2}+\|v(t)\|_\infty+\|Da(t)\|_\infty+
       \|g(t)\|_\infty\bigr)\le C_T .
\]
\end{proposition}

\begin{proof}
First work on $[0,S]$ with $S\le T$ and $S<S_n$. Let $B$ be a common
upper bound for $\|v\|_\infty$ given by Proposition~\ref{prop:poly-velocity}
or Proposition~\ref{prop:iso-velocity}. By Lemma~\ref{lem:log-lower} and
$h(0)-H\ge-2L$,
\begin{equation*}
 a(t,x)=\ee^{h(t,x)-H(x)}
 \ge\exp\bigl(-2L-4(L+B)\sqrt{T/\pi}-TB^2/4\bigr)=:c_T .
\end{equation*}
For the upper bound use the first component of \eqref{eq:mild}. With
$\mathcal F^a=(2d-v)a$ and
$\mathcal J^a=(|d|^2-\Delta H-d\cdot v)a$,
\begin{equation*}
 \|a(t)\|_\infty\le\ee^{2L}+
 \int_0^t \bigl(K_1(t-s)^{-1/2}+K_0\bigr)\|a(s)\|_\infty\,ds ,
\end{equation*}
where $K_1=2(2L+B)/\sqrt\pi$ and
$K_0=L^2+\|\Delta H\|_\infty+LB$. Choose $\nu>0$ so that
$K_1\sqrt{\pi/\nu}+K_0/\nu\le1/2$. Multiplication by $\ee^{-\nu t}$
and a supremum over $[0,S]$ give
$\sup_{[0,S]}\ee^{-\nu t}\|a(t)\|_\infty\le2\ee^{2L}$.
Since $b=av$, this gives $a\ge c_T$ and $\|U\|_\infty\le A_T$ on
$[0,S]\times\R^3$, with $A_T$ independent of $S$ and $n$.
Proposition~\ref{prop:local} now supplies a restart time
$\tau_0(A_T,c_T)>0$ independent of $S$ and $n$. If $S_n\le T$,
restarting at any $S<S_n$ with $S_n-S<\tau_0$ would extend the solution
beyond $S_n$. Thus $S_n>T$ for all $n$, and the estimates hold on
$[0,T]$, proving \eqref{eq:U-bounds}.

We next estimate $Da$ on this common interval. By \eqref{eq:flux},
\eqref{eq:source}, \eqref{eq:weighted-field} and \eqref{eq:U-bounds},
all components of $\mathcal F_j$ and $\mathcal J$ are uniformly bounded.
For $k=0,1$ and $0<\theta<1$, Gaussian scaling and the fundamental
theorem of calculus give
\[
 \|D^kG_r(\cdot+y)-D^kG_r\|_1
 \le \min\{2\|D^kG_r\|_1,|y|\|D^{k+1}G_r\|_1\}
 \lesssim |y|^\theta r^{-(k+\theta)/2}.
\]
Indeed, $\|D^kG_r\|_1\lesssim r^{-k/2}$ by scaling, and
$\min\{1,|y|r^{-1/2}\}\le|y|^\theta r^{-\theta/2}$ for $0<\theta<1$.
Fix $\theta=1/2$. By Lemma~\ref{lem:approx-data}(\ref{itm:normalized-data}),
$[U(0)]_\theta$ is uniformly bounded. Applying the translation bounds
to \eqref{eq:mild} yields
\[
 [U(t)]_\theta\le [U(0)]_\theta+
 C_T\int_0^t \bigl[(t-s)^{-(1+\theta)/2}+(t-s)^{-\theta/2}\bigr]\,ds
 \le C_T .
\]
In particular $\mathcal F^a=2da-b$ is uniformly $\theta$-H\"older
since $d$ is bounded and Lipschitz. Each entry of $D^2G_r$ has zero
integral, so
\begin{equation}\label{eq:cancel-gauss}
 \|D^2G_r*\mathcal F^a\|_\infty
 \le [\mathcal F^a]_\theta\int |D^2G_r(y)|\,|y|^\theta\,dy
 \le C_T r^{-1+\theta/2}.
\end{equation}
Differentiating the mild equation for $a$ now gives
\[
 \|Da(t)\|_\infty\le\|Da(0)\|_\infty+
 C_T\int_0^t\bigl[(t-s)^{-1+\theta/2}+(t-s)^{-1/2}\bigr]\,ds\le C_T .
\]
Indeed, the differentiated integrals truncated at $t-\varepsilon$
converge uniformly in space by these integrable bounds, and
$Da(0)=a_{0,n}(g_{0,n}-d)$ is uniformly bounded. Hence
$g=d+Da/a$ is uniformly bounded as well.
\end{proof}

The uniform bounds in \eqref{eq:U-bounds} allow us to iterate the
stability estimate \eqref{eq:stability} with a time step independent
of $n$. This compares the normalized solutions throughout $[0,T]$.
To identify $g=d+Da/a$ in the limit, we also need convergence of $Da_n$.

\begin{proposition}\label{prop:low-limit}
For every $T>0$, the sequence $U_n=(a_n,b_n)$ converges in $C_t\BUC$
to the mild solution $U=(a,b)$ of \eqref{eq:mild} with initial data
$(a_0,b_0)$, and $Da_n\to Da$ in $C_t\BUC$. Set
$\rho=\omega a$, $h=H+\log a$, $v=b/a$, $g=Dh=d+Da/a$ and $u=v-g$.
Then
\begin{equation}\label{eq:strong-low}
 v_n\longrightarrow v,\quad g_n\longrightarrow g
 \quad\text{in }C_t\BUC,\qquad
 \rho_n\longrightarrow\rho\quad\text{in }C_t(L^1\cap L^\infty).
\end{equation}
For every $t\in[0,T]$, $\int\rho(t)\,dx=M$, and
\begin{equation}\label{eq:limit-bounds}
 c_T\omega\le\rho(t)\le A_T\omega,\qquad
 \|v(t)\|_\infty+\|g(t)\|_\infty\le C_T,\qquad 0\le t\le T .
\end{equation}
The limit $U$ solves \eqref{eq:U-system} in distributions, and
$\rho$ and $m=\rho v$ satisfy \eqref{eq:heat-rho} and \eqref{eq:heat-m}
in distributions.
\end{proposition}

\begin{proof}
Since all approximations solve the same mild equation with the same
weight, the stability estimate \eqref{eq:stability} can be iterated on a
partition of $[0,T]$ whose mesh depends only on $A_T,c_T$ and the fixed
parameters. It gives
\[
 \sup_{0\le t\le T}\|U_n(t)-U_m(t)\|_\infty
 \le C_T\|U_{0,n}-U_{0,m}\|_\infty .
\]
Lemma~\ref{lem:approx-data}(\ref{itm:normalized-data}) therefore gives convergence to
$U=(a,b)$ in $C([0,T];\BUC)$. To pass to $Da_n$, observe that
$\mathcal F_n^a=2da_n-b_n$ converges uniformly. By
Proposition~\ref{prop:relative-continuation} and the bounds on $d$, these
fluxes are uniformly bounded in $C^\theta$ with $\theta=1/2$.
For $0<\beta<\theta$, the inequality
$[f]_\beta\le(2\|f\|_\infty)^{1-\beta/\theta}
[f]_\theta^{\beta/\theta}$ thus gives convergence of the flux in
$C^\beta$, uniformly in time. The sources $\mathcal J_n^a$ converge
uniformly as well. Using \eqref{eq:cancel-gauss} with exponent $\beta$
in the differentiated equations for $a_n-a_m$ gives
\[
 \sup_{[0,T]}\|Da_n-Da_m\|_\infty
 \le\|Da_{0,n}-Da_{0,m}\|_\infty
   +C_T\sup_{[0,T]}\bigl([\mathcal F_n^a-\mathcal F_m^a]_\beta+
                       \|\mathcal J_n^a-\mathcal J_m^a\|_\infty\bigr).
\]
The initial term tends to zero by
$Da_{0,n}=a_{0,n}(g_{0,n}-d)$ and Lemma~\ref{lem:approx-data}.
The uniform limit of $Da_n$ is $Da$, as follows by testing locally
against smooth functions.

The bound $a_n\ge c_T$ from Proposition~\ref{prop:relative-continuation}
allows passage to the limit in $v_n=b_n/a_n$ and $g_n=d+Da_n/a_n$.
Multiplication by $\omega\in L^1\cap L^\infty$ gives the density
convergence in \eqref{eq:strong-low}. In particular,
$\int\rho(t)\,dx=\lim_n M_n=M$, and the uniform estimates of
Proposition~\ref{prop:relative-continuation} give \eqref{eq:limit-bounds}.
On the common range \eqref{eq:U-bounds}, the local fluxes and sources
are uniformly Lipschitz in $U$. The field product satisfies the same
estimate by \eqref{eq:field-lipschitz}. All terms in \eqref{eq:mild}
therefore pass to the limit, since $\|DG_r\|_1\lesssim r^{-1/2}$.
Thus $U$ is the mild solution with initial data $(a_0,b_0)$ and solves
\eqref{eq:U-system} in distributions. Multiplication by the smooth weight
$\omega$ and the product identities in Lemma~\ref{lem:normalized-system}
give \eqref{eq:heat-rho} and \eqref{eq:heat-m} in distributions.
\end{proof}

\section{Higher-order estimates and strong traces}\label{sec:higher}

Fix $T>0$ and assume the hypotheses on the initial data in
Theorem~\ref{thm:main}. Proposition~\ref{prop:relative-continuation} gives uniform
supremum bounds for $v_n$ and $g_n$ on $[0,T]$. To estimate their
derivatives in unweighted Sobolev spaces, we subtract the smooth parts
of the initial profiles and the isothermal pressure tail. With the same
profiles for every approximation, set
\begin{equation*}
 d=\eta*g_0,\qquad \bar v(t)=\eta*v_0-\sigma t\,d,\qquad
 w_n=v_n-\bar v,\qquad z_n=g_n-d .
\end{equation*}
For the limit fields in \eqref{eq:strong-low}, write
$w=v-\bar v$ and $z=g-d$. In calculations for a single approximation
we suppress $n$.

For every $k\ge1$, $D^k\bar v$ and $D^kd$ are bounded in $C_tL^2$,
since one derivative can be placed on $v_0$ or $g_0$ and the remaining
derivatives on $\eta$. The profiles and all their spatial derivatives
are bounded in $C_t\BUC$.

\subsection{The shifted heat system}\label{subsec:shifted-system}

Put $c^2=\gamma\rho^{\gamma-1}$, so $c^2=1$ when $\gamma=1$.
Subtracting $\bar v$ in the effective momentum equation gives
\begin{equation}\label{eq:w-pde}
 (\partial_t-\Delta)w
 =\Delta\bar v-(v-2g)\cdot Dv-(c^2g-\sigma d)-E=:F_1 .
\end{equation}
Taking the gradient of \eqref{eq:heat-h} gives
\begin{equation}\label{eq:g-derivation}
 (\partial_t-\Delta)g+D\Div v=(Dg)^T(2g-v)-(Dv)^Tg ,
\end{equation}
and therefore
\begin{equation}\label{eq:z-pde}
 (\partial_t-\Delta)z+D\Div w
 =\Delta d-D\Div\bar v+(Dg)^T(2g-v)-(Dv)^Tg=:F_2 .
\end{equation}

The energy estimate for \eqref{eq:w-pde} and \eqref{eq:z-pde} involves
the $H^2$ norm at the endpoints of a time interval. We first identify a
continuous representative for which these endpoint values are defined.

\begin{lemma}\label{lem:strong-time-trace}
If $y\in L^2_tH^3$ and $\partial_t y\in L^2_tH^1$, then $y$ has a
representative in $C_tH^2$.
\end{lemma}

\begin{proof}
Let $P=\Pi_{N'}-\Pi_N$ for $N'>N$. The band-limited function $Py$ belongs
to $H^1(0,T;H^2)$. Integrate the derivative of $\|Py\|_{H^2}^2$
starting at a time when it is at most its mean. The Fourier pairing
between $H^3$ and $H^1$ gives
\begin{equation}\label{eq:trace}
 \sup_{0\le t\le T}\|Py(t)\|_{H^2}^2
 \le T^{-1}\|\Pi_{>N}y\|_{L^2_tH^2}^2
 +2\|\Pi_{>N}y\|_{L^2_tH^3}\|\Pi_{>N}\partial_ty\|_{L^2_tH^1}.
\end{equation}
The right-hand side tends to zero uniformly in $N'>N$. Hence $\Pi_Ny$
converges in $C_tH^2$ to a representative of $y$.
\end{proof}

In the shifted system, the second equation contains the second-order
coupling $D\Div w$. To estimate this term along with the diffusion in
both equations, we use the energy and dissipation of a pair
$Y=(y_1,y_2)$, writing
$\mathcal A(Y)=\sum_i\|y_i\|_{H^2}^2$ and
$\mathcal D(Y)=\sum_i\|Dy_i\|_{H^2}^2$.

\begin{lemma}\label{lem:coupled-heat-energy}
Let $Y=(y_1,y_2)\in L^2_tH^3$ be a distributional solution of
\begin{equation*}
 (\partial_t-\Delta)y_1=f_1,\qquad
 (\partial_t-\Delta)y_2+D\Div y_1=f_2 .
\end{equation*}
If $f_1,f_2\in L^2_tH^1$, then $Y\in C_tH^2$
and, for almost every $t\in(0,T)$,
\begin{equation}\label{eq:linear-differential}
 \frac{d}{dt}\mathcal A(Y)+\tfrac12\mathcal D(Y)
 \le \tfrac12\mathcal A(Y)+2\sum_i\|f_i\|_{H^1}^2 .
\end{equation}
\end{lemma}

\begin{proof}
The equations give $\partial_tY\in L^2_tH^1$, so
Lemma~\ref{lem:strong-time-trace} gives $Y\in C_tH^2$.
Testing the two equations in $H^2$, we bound the coupling term
$-\langle D\Div y_1,y_2\rangle_{H^2}$ in modulus
by $\|Dy_1\|_{H^2}\|Dy_2\|_{H^2}\le\mathcal D(Y)/2$.
Each source pairing is bounded by
$\|f_i\|_{H^1}\|y_i\|_{H^3}
 \le \|f_i\|_{H^1}^2+\tfrac14(\|y_i\|_{H^2}^2+\|Dy_i\|_{H^2}^2)$.
Fourier truncation justifies these pairings and the time differentiation.
In the integrated identity the endpoint terms converge by \eqref{eq:trace},
and the time integrals converge by Cauchy--Schwarz in $L^2_t$.
\end{proof}

\subsection{Spatial integrability of the shifted variables}\label{subsec:identify}

On $\R^3$, bounded spatial derivatives alone give no control of their
$L^2$ norms at infinity. Before using Lemma~\ref{lem:coupled-heat-energy}
for $(w_n,z_n)$, we justify the required spatial integrability for each
fixed $n$.

\begin{lemma}\label{lem:shifted-integrability}
For each smooth approximation on $[0,T]$ from Proposition~\ref{prop:relative-continuation},
its shifted pair $Y_n=(w_n,z_n)$ belongs to $C_tH^2\cap L^2_tH^3$,
with $\partial_tY_n\in L^2_tH^1$. Its $H^2$ initial trace is
$(v_{0,n}-\eta*v_0,g_{0,n}-d)$.
\end{lemma}

\begin{proof}
Fix $n$ and write $Y=(w,z)$ for the shifted pair.
Expand $v=w+\bar v$ and $g=z+d$ in the differentiated factors in
\eqref{eq:w-pde}--\eqref{eq:z-pde}, keeping the coefficients at their
values along this approximation. The resulting system has the form
\begin{equation}\label{eq:frozen-system}
 \partial_tY-\mathcal P(D)Y=A_j\partial_jY+A_0Y+f,\qquad
 \mathcal P(D)Y=(\Delta w,\Delta z-D\Div w),
\end{equation}
where the coefficients $A_j,A_0$ are bounded in $W^{1,\infty}$
uniformly on $[0,T]$ for fixed $n$. In fact the first-order coefficients
are linear combinations of components of $v,g$, and the zeroth-order
term is $(-c^2z,0)$. We have $Dc^2=(\gamma-1)c^2g$.
The prescribed source is
\begin{align*}
 f_1={}&\Delta\bar v-(v-2g)\cdot D\bar v-(c^2-\sigma)d-E,\\
 f_2={}&\Delta d-D\Div\bar v+(Dd)^T(2g-v)-(D\bar v)^Tg .
\end{align*}
Here and below the $W^{1,\infty}$ bounds on the fixed approximation
follow from Corollary~\ref{cor:approx-smooth}.

We check that $f\in L^2_tH^1$. The profile derivatives have the
required spatial integrability, and
$\|E\|_{H^1}\lesssim_\ell\|\rho\|_1+\|\rho\|_2$ follows from the
multiplier estimate \eqref{eq:field-sobolev} with $k=1$. The remaining term
$(c^2-\sigma)d$ vanishes at $\gamma=1$. At $\gamma>1$, write
\begin{equation}\label{eq:background-pressure}
 c^2d=c^2(d-g_{0,n})+(a/a_{0,n})^{\gamma-1}p_{0,n},
 \qquad p_{0,n}=\gamma\rho_{0,n}^{\gamma-1}g_{0,n}.
\end{equation}
Both summands belong to $L^\infty_tL^2$ by Lemma~\ref{lem:approx-data}
and \eqref{eq:U-bounds}. The identity
$\partial_j(c^2d)=c^2\partial_jd+(\gamma-1)g_jc^2d$ gives the same
bound for the first spatial derivative.

For $R\ge1$, let $\vartheta_R(x)=\exp(-\sqrt{1+|x|^2}/R)$ and
$Y_R=\vartheta_RY$. For this fixed $n$, $Y$ has bounded spatial
derivatives of all orders on $[0,T]$. Hence $Y_R\in C_tH^2\cap L^2_tH^3$.
Multiplication of \eqref{eq:frozen-system} by $\vartheta_R$ replaces
each $\partial_j$ acting on $Y$ by $\partial_j-\partial_j\log\vartheta_R$
acting on $Y_R$. Consequently the principal operator acting on $Y_R$
is still $\mathcal P(D)$, and the additional terms have order at most one.
The derivatives of $\log\vartheta_R$ of orders one through three are
bounded uniformly for $R\ge1$. Thus the coefficients of all lower-order
terms and their first spatial derivatives are bounded by $C_{n,T}$
independently of $R$. Also
$\|\vartheta_R f\|_{H^1}\le C\|f\|_{H^1}$. Applying
\eqref{eq:linear-differential} to $Y_R$ gives
\[
 \frac{d}{dt}\mathcal A(Y_R)+\tfrac12\mathcal D(Y_R)
 \le C_{n,T}\mathcal A(Y_R)+C\|f\|_{H^1}^2 .
\]
Since $\|\vartheta_R\|_{W^{2,\infty}}$ is bounded uniformly and
$Y(0)\in H^2$ by Lemma~\ref{lem:approx-data}(\ref{itm:shifted-data}), Gr\"onwall's inequality
bounds $\|Y_R\|_{L^\infty_tH^2}+\|Y_R\|_{L^2_tH^3}$ independently of $R$.
As $R\to\infty$, the weights tend to one smoothly on compact sets.
Weak compactness in $L^2_tH^3$ and local convergence identify the limit
as $Y$. Thus $Y\in L^2_tH^3$, and \eqref{eq:frozen-system} gives
$\partial_tY\in L^2_tH^1$. Lemma~\ref{lem:strong-time-trace} gives
$Y\in C_tH^2$. This representative agrees at every time with the
original smooth pair by local continuity, which identifies the initial trace.
\end{proof}

\subsection{Estimates independent of the approximation}\label{subsec:uniform}

The preceding argument allows constants depending on $n$. For an
estimate uniform in $n$, we use the supremum bounds on $v_n$ and $g_n$
to control the products of first derivatives that arise when
differentiating the transport terms.

\begin{lemma}\label{lem:gradient-interpolation}
For a bounded vector field $q$ with $Dq\in H^1$, we have
\begin{equation}\label{eq:grad-four-squared}
 \|Dq\|_4^2\le (2+\sqrt3)\|q\|_\infty\|D^2q\|_2 .
\end{equation}
\end{lemma}

\begin{proof}
We have $Dq\in L^2\cap L^6\subset L^3\cap L^4$.
Integrate by parts in
$\int\chi_R\,\partial_jq_i\,\partial_jq_i\,|Dq|^2\,dx$,
placing one derivative on the other factors. The cutoff error is bounded
by $CR^{-1}\|q\|_\infty\|Dq\|_3^3$. Since
$|\Delta q|\le\sqrt3|D^2q|$ and $|D(|Dq|^2)|\le2|Dq||D^2q|$,
the remaining integrals are bounded by
$(2+\sqrt3)\|q\|_\infty\|Dq\|_4^2\|D^2q\|_2$.
Letting $R\to\infty$ proves \eqref{eq:grad-four-squared}, first for
smooth $q$ and then by spatial mollification.
\end{proof}

Lemma~\ref{lem:gradient-interpolation} controls the transport products
in terms of the $H^2$ energy of the shifted variables. For the pressure
term, the decomposition \eqref{eq:background-pressure} relates the
part involving $d$ to the initial data.

\begin{proposition}\label{prop:uniform-higher}
For each smooth approximation, the sources in
\eqref{eq:w-pde} and \eqref{eq:z-pde} satisfy
\begin{equation}\label{eq:source-bound}
 \|F_{1,n}(t)\|_{H^1}+\|F_{2,n}(t)\|_{H^1}
 \le C_T\bigl(1+\|w_n(t)\|_{H^2}+\|z_n(t)\|_{H^2}\bigr).
\end{equation}
The shifted variables satisfy
\begin{equation}\label{eq:uniform-higher}
 \|(w_n,z_n)\|_{C_tH^2}+\|(w_n,z_n)\|_{L^2_tH^3}
 +\|(\partial_tw_n,\partial_tz_n)\|_{L^2_tH^1}\le C_T ,
\end{equation}
where $C_T$ is independent of $n$.
\end{proposition}

\begin{proof}
Fix an approximation and suppress $n$. By
Lemma~\ref{lem:shifted-integrability}, the shifted variables have the
regularity required to apply Lemma~\ref{lem:coupled-heat-energy}.
For $q_1,q_2\in\{v,g\}$, the uniform supremum bounds and
Lemma~\ref{lem:gradient-interpolation} give
\begin{align*}
 \|q_1\otimes Dq_2\|_{H^1}
 &\le C\bigl(\|q_1\|_\infty\|Dq_2\|_{H^1}
                   +\|Dq_1\|_4\|Dq_2\|_4\bigr)\\
 &\le C_T\bigl(1+\|w\|_{H^2}+\|z\|_{H^2}\bigr).
\end{align*}
In the second inequality we used $2ab\le a^2+b^2$, together with the
$H^1$ bounds for $D\bar v$ and $Dd$. To estimate the pressure at
$\gamma>1$, use $c^2g=c^2z+c^2d$ and \eqref{eq:background-pressure}.
The bounds for $d-g_{0,n}$ and $p_{0,n}$ in $L^2$, for
$a/a_{0,n}$ in $L^\infty$, and for $c^2$ and $g$ in $L^\infty$
are all uniform in $n$. The derivative formula following
\eqref{eq:background-pressure} therefore also bounds
$\|c^2d\|_{H^1}$ uniformly. Since
$\|Dc^2\|_\infty\le(\gamma-1)\|c^2\|_\infty\|g\|_\infty$,
\[
 \|c^2g\|_{H^1}\le C_T(1+\|z\|_{H^1}),\qquad \gamma>1 .
\]
At $\gamma=1$ the shifted pressure term is exactly $c^2g-\sigma d=z$.
Finally, $\|E\|_{H^1}\lesssim_\ell\|\rho\|_1+\|\rho\|_2$ is bounded
uniformly because $\rho\le A_T\omega$. Inserting these bounds into
\eqref{eq:w-pde} and \eqref{eq:z-pde} proves \eqref{eq:source-bound}.
Applying \eqref{eq:linear-differential} to $Y=(w,z)$ now gives
$\frac{d}{dt}\mathcal A(Y)+\frac12\mathcal D(Y)
 \le C_T(1+\mathcal A(Y))$. Gr\"onwall's inequality and the uniform
$H^2$ bound on $Y(0)$ from Lemma~\ref{lem:approx-data}(\ref{itm:shifted-data}) control
$\sup_t\mathcal A(Y)$ and $\int_0^T\mathcal D(Y)\,dt$.
The equations and \eqref{eq:source-bound} then control
$\partial_tY$ in $L^2_tH^1$, proving \eqref{eq:uniform-higher}.
\end{proof}

\subsection{Passage to the limit and strong traces}\label{subsec:limit}

The profiles $\bar v$ and $d$ are fixed, so the convergence in
\eqref{eq:strong-low} identifies the weak limits furnished by
\eqref{eq:uniform-higher} with $v-\bar v$ and $g-d$. It remains to
pass to the shifted equations and identify the initial traces.

\begin{proposition}\label{prop:shifted-limit}
Let $v,g$ be the limits in \eqref{eq:strong-low}. The shifted pair
$(w,z)=(v-\bar v,g-d)$ solves \eqref{eq:w-pde}--\eqref{eq:z-pde}
in distributions and satisfies \eqref{eq:uniform-higher} with
$(w_n,z_n)$ replaced by $(w,z)$. Its strong $H^2$ traces at $t=0$ are
\begin{equation*}
 w(0)=v_0-\eta*v_0,\qquad z(0)=g_0-d .
\end{equation*}
\end{proposition}

\begin{proof}
By \eqref{eq:uniform-higher}, a subsequence of $(w_n,z_n)$ converges
weakly in $L^2_tH^3$ and weakly-* in $L^\infty_tH^2$. Its time
derivatives converge weakly in $L^2_tH^1$ along a further subsequence.
The uniform convergence in \eqref{eq:strong-low} identifies the limits
as $(w,z)=(v-\bar v,g-d)$ and their distributional time derivatives.

The pair $(w,z)$ satisfies \eqref{eq:w-pde}--\eqref{eq:z-pde} in
distributions. Indeed, in each transport term, the undifferentiated
factor converges uniformly, while the derivative factor converges
weakly in $L^2$ on the support of any test function. The normalized
pressure gradients $c_n^2g_n$ converge uniformly because $\rho_n,g_n$
converge uniformly and the densities are bounded. The bound \eqref{eq:field-sobolev} with $k=1$
and the convergence in \eqref{eq:strong-low} give $E_n\to E$ in $C_tH^1$,
since
\[
 \|E_n-E\|_{H^1}\lesssim_\ell
      \|\rho_n-\rho\|_1+\|\rho_n-\rho\|_2
\]
and $\|f\|_2^2\le\|f\|_1\|f\|_\infty$.
The remaining terms are linear, with profiles independent of $n$.

Lemma~\ref{lem:strong-time-trace} gives strong time continuity of both
shifted variables in $H^2$. Weak lower semicontinuity in the bounds of
Proposition~\ref{prop:uniform-higher} gives the asserted estimate,
since the essential supremum of the $H^2$ norm of a continuous path is
its supremum. The $C_tH^2$ representatives also belong to $C_t\BUC$
and hence coincide at every time with the uniform limits in
\eqref{eq:strong-low}. Evaluation at $t=0$ and
Lemma~\ref{lem:approx-data}(\ref{itm:shifted-data}) give the stated initial traces.
\end{proof}

\section{Proof of the main theorem}\label{sec:proof}

\begin{proof}[Proof of Theorem~\ref{thm:main}]
Fix $T>0$. Let $U=(a,b)$ be the limit from
Proposition~\ref{prop:low-limit}, and set
\begin{equation*}
 \rho=\omega a,\qquad v=b/a,\qquad g=D\log\rho=d+Da/a,\qquad
 u=v-g,\qquad E=\mathcal N[\rho].
\end{equation*}
By \eqref{eq:strong-low} and \eqref{eq:limit-bounds},
$\rho\in C_t(L^1\cap L^\infty)$, $\int\rho(t)\,dx=\int\rho_0\,dx$,
and $c_T\omega\le\rho\le A_T\omega$. Comparison of $\omega$ with $\rho_0$
gives \eqref{eq:reg-relative}. The convergence at $t=0$ gives
$a(0)=a_0$ and $\rho(0)=\rho_0$.

\smallskip
\noindent\emph{Velocity and logarithmic gradient.}
Write $w=v-\bar v$ and $z=g-d$, where
$\bar v=\eta*v_0-\sigma td$. Proposition~\ref{prop:shifted-limit} gives
$w,z\in C_tH^2\cap L^2_tH^3$ with
$w(0)=v_0-\eta*v_0$ and $z(0)=g_0-d$. In particular $v(0)=v_0$,
$g(0)=g_0$ and $u(0)=u_0$. The identities
\begin{align}
 v-v_0+\sigma t g_0&=w-w(0)+\sigma t z(0),&
 g-g_0&=z-z(0),\label{eq:undo-v}\\
 u+\sigma td&=w-z+\eta*u_0 \label{eq:u-decomposition}
\end{align}
give \eqref{eq:reg-shifts}, since $g_0-d=z(0)\in H^2$.
The profiles $d$ and $\bar v$ have bounded spatial derivatives of
every order. Their derivatives of positive order are in $L^2$, uniformly
on $[0,T]$, by convolution and $Dg_0,Du_0\in H^1$. Thus
\begin{equation}\label{eq:g-regularity}
 Du,\ Dv,\ Dg\in C_tH^1\cap L^2_tH^2,\qquad
 u,\ v,\ g\in C_t\BUC .
\end{equation}
Also $\eta*u_0\in H^3$, so \eqref{eq:u-decomposition} yields
$u+\sigma td\in L^2_tH^3$, and hence $u\in L^2_tH^3$ for $\gamma>1$.

\smallskip
\noindent\emph{Density and mass flux.}
We use $D\rho=\rho g$ to recover the density derivatives. First
$\rho\in C_tL^2$ by interpolation between $L^1$ and $L^\infty$.
Since $g\in C_tL^\infty$, we obtain $D\rho\in C_tL^2$. The weak
product rule $D(\rho g)=g\otimes D\rho+\rho Dg$ then gives
$\rho\in C_tH^2$. For $k=2,3$, $H^k(\R^3)$ is a multiplication
algebra, and multiplication by $d\in W^{k,\infty}$ acts boundedly on
$H^k$. Since $g=d+z$, it follows successively that
\[
 \rho g\in C_tH^2,\qquad \rho\in C_tH^3,\qquad
 \rho g\in L^2_tH^3,\qquad \rho\in L^2_tH^4 .
\]
In the third assertion we used $z\in L^2_tH^3$ and
$\sup_{t\le T}\|\rho(t)\|_{H^3}<\infty$. Applying the same product
bounds to \eqref{eq:u-decomposition} gives
$\rho u\in C_tH^2\cap L^2_tH^3$. Substitution of
$\rho v=\rho u+D\rho$ in \eqref{eq:heat-rho} therefore gives
\begin{equation}\label{eq:rho-t}
 \rho_t=-\Div(\rho u)\in C_tH^1\cap L^2_tH^2 .
\end{equation}
In particular the $H^3$ trace of $\rho$ at zero is $\rho_0$.

\smallskip
\noindent\emph{Electric field.}
Apply \eqref{eq:field-sobolev} with $k=4$ to $\rho(t)-\rho(s)$ and
with $k=5$ to $\rho(t)$. The density estimates give
$E\in C_tH^4\cap L^2_tH^5$. In particular $E\in C_tW^{1,\infty}$.
Define
\[
 V=(4\pi\ell^2|\cdot|)^{-1}*\rho .
\]
Splitting at $|\xi|=1$ and applying Cauchy--Schwarz give
$\||\xi|^{-2}\widehat\varphi\|_1\le
 C(\|\varphi\|_1+\|\varphi\|_2)$, since $|\xi|^{-2}\in
 L^2(\{|\xi|>1\})$. Hence $\widehat V\in C_tL^1_\xi$, and Fourier
inversion yields $V\in C_tC_0$. We have $DV=E$ and
$-\ell^2\Delta V=\rho$. The harmonic difference of any two potentials
with this normalization vanishes by the maximum principle.

\smallskip
\noindent\emph{Time derivatives.}
Put $p=\gamma\rho^{\gamma-1}g$. For $\gamma>1$,
\begin{equation*}
 p=\gamma\rho^{\gamma-1}(g-g_0)
     + (a/a_0)^{\gamma-1}p_0,\qquad
 Dp=\gamma\rho^{\gamma-1}Dg+(\gamma-1)p\otimes g .
\end{equation*}
Both coefficients in the first identity are continuous in time in the
supremum norm, and $p_0\in L^2$. It follows first that $p\in C_tL^2$
and then, by the second identity, that $p\in C_tH^1$. At $\gamma=1$,
$p-g_0=g-g_0\in C_tH^2$. Thus $p-\sigma g_0\in C_tH^1$ for every
exponent in the theorem. Equation \eqref{eq:effective} for $v$ and
the gradient equation \eqref{eq:g-derivation} now read
\begin{align}
 v_t+\sigma g_0
   &=\Delta v-(v-2g)\cdot Dv-(p-\sigma g_0)-E,
       \label{eq:time-recover-v}\\
 g_t&=\Delta g-D\Div v+(Dg)^T(2g-v)-(Dv)^Tg .
       \label{eq:time-recover-g}
\end{align}
All terms on the right are in $C_tL^2\cap L^2_tH^1$. Indeed, the
linear second derivatives have this regularity by \eqref{eq:g-regularity}.
Each transport product has a bounded continuous factor multiplying a
first derivative in $C_tL^2$. Its spatial derivatives are sums of a
bounded factor times a second derivative and products of two first
derivatives in $C_tL^4$, since $H^1(\R^3)\hookrightarrow L^4(\R^3)$.
Subtracting the two identities gives the same regularity for
$u_t+\sigma g_0$, as asserted in \eqref{eq:reg-time}. Integration
of \eqref{eq:rho-t} and \eqref{eq:time-recover-v}--\eqref{eq:time-recover-g}
in their respective Hilbert spaces shows that $\rho\in C^1([0,T];H^1)$
and that each shifted variable in \eqref{eq:reg-shifts} belongs to
$C^1([0,T];L^2)$. For example,
\[
 \rho_t(0)=-\Div(\rho_0u_0),\qquad
 g_t(0)=-D\Div u_0-D(u_0\cdot g_0),
\]
where $\Delta g_0=D\Div g_0$ follows from $g_0=D\log\rho_0$.

\smallskip
\noindent\emph{Momentum equation.}
Equation \eqref{eq:heat-m} holds for the limit in distributions.
Using $\rho(v-u)=D\rho$, we may write it as
\[
 \partial_t(\rho v_i)+\partial_j(\rho u_jv_i)
 =\partial_j(\rho\,\partial_jv_i)-\partial_i\rho^\gamma-\rho E_i .
\]
Differentiating the mass equation and using $\rho g_i=\partial_i\rho$ gives
\[
 \partial_t(\rho g_i)+\partial_j(\rho u_jg_i)
 =-\partial_j(\rho\,\partial_i u_j).
\]
Subtraction and $v=u+g$ yield
\begin{equation*}
 \partial_t(\rho u_i)+\partial_j(\rho u_i u_j)+\partial_i\rho^\gamma
 =\partial_j\!\left[\rho(\partial_j u_i+\partial_i u_j)
                    +\rho\,\partial_jg_i\right]-\rho E_i .
\end{equation*}
Here $Dg=D^2\log\rho$ on each compact cylinder, where
$\rho\ge c_T\omega$ has a strictly positive minimum. The regularity
just proved makes all derivatives in \eqref{eq:model} locally
integrable, so the distributional identities hold almost everywhere.

\smallskip
\noindent\emph{Uniqueness in the stated class.}
Let $(\tilde\rho,\tilde u,\tilde V)$ be a positive strong solution with
the same initial data and regularity \eqref{eq:regularity}. Write
$\tilde g=D\log\tilde\rho$ and $\tilde v=\tilde u+\tilde g$.
On every compact cylinder, $\tilde\rho$ is continuous and bounded
away from zero. Its mass equation and initial trace imply
\begin{equation*}
 \log\frac{\tilde\rho(t)}{\rho_0}
 =-\int_0^t(\tilde u\cdot\tilde g+\Div\tilde u)(s)\,ds .
\end{equation*}
The integral on the right belongs to $C_t\BUC$ globally, since
$\tilde u,\tilde g\in C_t\BUC$ and
$\Div\tilde u\in L^2_tH^2\subset L^1_t\BUC$. Consequently
$\tilde a=\tilde\rho/\omega=a_0\exp(\log(\tilde\rho/\rho_0))$
and $\tilde b=\tilde a\tilde v$ belong to $C_t\BUC$, with
$\inf_{[0,T]\times\R^3}\tilde a>0$. This also gives
$\tilde\rho\in C_tL^1$. By the uniqueness of the normalized Poisson
potential, $D\tilde V=\mathcal N[\omega\tilde a]$. The initial density
trace and local positivity identify the trace of $\tilde g$ with
$g_0$. Hence $(\tilde a,\tilde b)(0)=(a_0,b_0)$.

The product identities from Lemmas~\ref{lem:effective-system}
and~\ref{lem:normalized-system} apply locally
to the strong solution and give \eqref{eq:U-system} for
$(\tilde a,\tilde b)$ in distributions. The bounded fluxes and source, the $\BUC$ time trace, and the
heat-kernel representation in Lemma~\ref{lem:heat-duhamel} give
\eqref{eq:mild} for $(\tilde a,\tilde b)$. Both pairs are therefore mild
solutions with the same weight and initial values. Their bounds and
positive lower bounds for $a,\tilde a$ permit
Proposition~\ref{prop:local} to be applied on finitely many common
subintervals of $[0,T]$. Thus $\tilde a=a$ and $\tilde b=b$.
Recovering $\rho,v,g,u$ and the normalized potential proves uniqueness.

\smallskip
\noindent\emph{Behavior at infinity.}
Lemma~\ref{lem:initial-buc} gives $g_0\in C_0$, whence also
$d,\eta*v_0\in C_0$. Since $H^2(\R^3)\hookrightarrow C_0(\R^3)$,
$w,z\in C_tC_0$ and therefore $v,g,u\in C_tC_0$.
We likewise have $\rho,E,V\in C_tC_0$. A continuous path on $[0,T]$
has compact image in $C_0$, and every compact subset of $C_0$
vanishes uniformly at infinity, as follows by taking a finite
$\varepsilon$-net in the supremum norm. This proves the asserted
uniform limits.

The estimates \eqref{eq:limit-bounds} and \eqref{eq:field-sobolev},
together with the bounds for $w,z$ in Proposition~\ref{prop:shifted-limit}, give
\eqref{eq:reg-sup} with the dependence stated in the theorem.
For $5/4\le\gamma\le3$, Lemma~\ref{lem:initial-buc} controls
$\|p_0\|_2$ by the norms in \eqref{eq:initial}. Uniqueness identifies the
solutions on overlapping finite time intervals, giving the global solution.
\end{proof}

\appendix
\section{Local theory for the normalized system}\label{app:local}

\subsection{Properties of the initial data}\label{subsec:initial-buc}

The construction of the approximate data in Section~\ref{subsec:approx}
uses uniform continuity of $g_0$ and bounds on the spatial increments
of $a_0$ and $b_0$. We derive these properties from \eqref{eq:initial}
and determine when the same hypotheses control the initial pressure
gradient in $L^2$.

\begin{lemma}\label{lem:initial-buc}
Assume \eqref{eq:initial}. Then $g_0$ has a representative in
$C_0(\R^3)\cap L^6(\R^3)$ and
\begin{equation*}
 \|g_0\|_6\le C\|Dg_0\|_2,\qquad
 [g_0]_{1/2}+[a_0]_1+[b_0]_{1/2}\le K_0 .
\end{equation*}
Here $a_0=\rho_0/\omega$, $b_0=a_0v_0$, and $K_0$ depends only on
$\|g_0\|_\infty$, $\|Dg_0\|_{H^1}$, $\|u_0\|_{H^2}$ and the fixed
mollifier. The same seminorm bounds hold for
$g_{0,n},a_{0,n},b_{0,n}$, uniformly in $n$. For $\gamma\ge5/4$,
\begin{equation}\label{eq:automatic-pressure}
 \|\gamma\rho_0^{\gamma-1}g_0\|_2^2
 \le\frac{\gamma^2}{2(\gamma-1)}
       \|\rho_0^{2(\gamma-1)}\|_2\,\|\Div g_0\|_2 .
\end{equation}
\end{lemma}

\begin{proof}
The Sobolev inequality gives $Dg_0,Du_0\in L^6(\R^3)$. By Morrey's
inequality on unit balls and the boundedness of $g_0$ and $u_0$, these
fields have globally $1/2$-H\"older representatives with
\[
 [g_0]_{1/2}+[u_0]_{1/2}
 \le C\bigl(\|g_0\|_\infty+\|Dg_0\|_{H^1}+\|u_0\|_{H^2}\bigr).
\]
See \cite[Theorem 9.9 and Corollary 9.14]{Brezis}. For pairs of points
at distance at least one the estimate follows directly from the
supremum norms. By \eqref{eq:weight-bounds},
$Da_0=a_0(g_0-d)$ is bounded. A bounded Lipschitz function has finite
$1/2$-H\"older seminorm, so $b_0=a_0(u_0+g_0)$ has the asserted bound.
Convolution preserves the bounds for $g_{0,n}$ and $u_{0,n}$, while
$|h_{0,n}-H|\le2L$ and $Da_{0,n}=a_{0,n}(g_{0,n}-d)$ give the
corresponding estimates for $a_{0,n}$ and $b_{0,n}$.

To identify the behavior of $g_0$ at infinity, let $m_R$ be its average
over $B_R(0)$. The scale-invariant Sobolev--Poincar\'e inequality gives
\[
 \|g_0-m_R\|_{L^6(B_R)}\le C\|Dg_0\|_{L^2(B_R)} .
\]
Since $|m_R|\le\|g_0\|_\infty$, some sequence $R_j\to\infty$ has
$m_{R_j}\to m\in\R^3$. Fatou's lemma yields
$\|g_0-m\|_6\le C\|Dg_0\|_2$. Uniform continuity then gives
$g_0(x)\to m$ as $|x|\to\infty$. Indeed, a sequence escaping to infinity
on which $|g_0-m|\ge\varepsilon$ would give disjoint balls of a common
positive radius on which $|g_0-m|\ge\varepsilon/2$. If $m\ne0$, the
directional derivative of $h_0=\log\rho_0$ along the ray
$r\mapsto rm/|m|$ tends to $|m|$. This contradicts the upper bound
$h_0\le\log\|\rho_0\|_\infty$. Thus $m=0$.

Finally, put $\beta=\gamma-1\ge1/4$ and $f=\rho_0^{2\beta}$.
Then $\|f\|_2^2\le\|\rho_0\|_\infty^{4\beta-1}\|\rho_0\|_1$.
Since $Df=2\beta f g_0$, integration by parts gives
\[
 2\beta\int\chi_R f|g_0|^2\,dx
 =-\int\chi_R f\,\Div g_0\,dx-\int f g_0\cdot D\chi_R\,dx .
\]
The absolute value of the last integral is bounded by
$\|f\|_{L^2(|x|\ge R)}\|g_0\|_6\|D\chi_R\|_3$. The final factor is
independent of $R$, so the cutoff error tends to zero. Fatou's lemma
on the left and Cauchy--Schwarz on the right prove
\eqref{eq:automatic-pressure}.
\end{proof}

For $1<\gamma<5/4$, Theorem~\ref{thm:main} imposes
\eqref{eq:initial-pressure} separately. This additional assumption
also concerns time regularity, since the effective momentum equation
relates $p_0$ to the initial trace of $v_t$.

\begin{remark}\label{rem:pressure-threshold}
The threshold $\gamma=5/4$ in Lemma~\ref{lem:initial-buc} is sharp for
the implication from \eqref{eq:initial} to \eqref{eq:initial-pressure}.
Fix $0<\beta=\gamma-1<1/4$. Since $3\beta<3/4$, we may choose
$m>3$ and $0<\lambda<1/4$ so that $m\beta<\lambda+1/2$.
Set $u_0=0$ and, for $r=|x|\ge2$, let
\[
 \rho_0(x)=r^{-m}\exp\bigl(\sin(r^\lambda)\bigr).
\]
Choose a smooth positive radial extension to $r<2$ that is constant
near the origin. For $j=1,2,3$ one has
$|D^j\log\rho_0|=O(r^{j(\lambda-1)})$ at infinity.
In particular $g_0$ is bounded, and $Dg_0,D^2g_0\in L^2$ since
$\lambda<1/4$. The density is integrable, and its derivatives through
order three are $O(r^{-m})$, so \eqref{eq:initial} holds.
On the set $|\cos(r^\lambda)|\ge1/2$ we have
$|g_0|\ge c r^{\lambda-1}$ for all sufficiently large $r$. It follows
that the squared $L^2$ norm of $\rho_0^\beta g_0$ dominates a positive
multiple of
\[
 \int_{\{r>R:\,|\cos(r^\lambda)|\ge1/2\}} r^{2\lambda-2m\beta}\,dr .
\]
This integral diverges. Indeed, after $s=r^\lambda$ the exponent of $s$
is $(2\lambda-2m\beta+1)/\lambda-1>-1$, and the integration set
contains an interval of fixed positive length in every full period.

For $\gamma>1$ the condition $p_0\in L^2$ is also necessary for the
time trace in \eqref{eq:reg-time} under the other regularity and
initial trace assumptions of the theorem. To see this, solve
\eqref{eq:effective} for the normalized pressure gradient to obtain
\[
 \gamma\rho^{\gamma-1}g=\Delta v-(v-2g)\cdot Dv-E-v_t .
\]
All terms on the right have strong $L^2$ traces at zero by
\eqref{eq:regularity}. On each compact set the left-hand side converges
uniformly to $p_0$, by $\rho(t)\to\rho_0$ in $H^3$ and local positivity
of $\rho_0$. Thus the $L^2$ trace equals $p_0$.
\end{remark}

\subsection{Local existence and stability}\label{subsec:proof-local}

\begin{proof}[Proof of Proposition~\ref{prop:local}]
Fix the weight $\omega$ and equip $\BUC(\R^3;\R^4)$ with the maximum
of the componentwise supremum norms.
For $0<c\le A$, write
\[
 \Omega_{A,c}=\{U=(a,b)\in\BUC:\ \|U\|_\infty\le A,\ \inf_x a(x)\ge c\}.
\]
The maps $\mathcal F_j,\mathcal J$ in \eqref{eq:flux}--\eqref{eq:source}
are bounded and Lipschitz on $\Omega_{A,c}$. For the local terms this
follows from the bounds on $d,\Delta H,\omega^{\gamma-1}$ and the
bounded derivatives of $(a,b)\mapsto b_ib_j/a$ and $a\mapsto a^\gamma$
on the indicated range. Lemma~\ref{lem:field-reg} treats the product
$a\,\mathcal N[\omega a]$. Denote by $K$ a common bound for the sum of
their norms and by $\Lambda$ a bound for the sum of their Lipschitz constants. These
constants depend only on $A,c$ and the fixed parameters in the proposition.

Put $r=\min(1,\inf a_{\mathrm{in}}/2)$ and consider the closed ball
of radius $r$ about $U^0(t)=G_t*U_{\mathrm{in}}$ in
$C([0,\tau];\BUC)$, with $U^0(0)=U_{\mathrm{in}}$. By positivity and
unit mass of $G_t$, each path in the ball takes its values in
$\Omega_{A,c}$ with $A=\|U_{\mathrm{in}}\|_\infty+1$ and
$c=\inf a_{\mathrm{in}}/2$. Let $\mathcal T U$ denote the right-hand
side of \eqref{eq:mild}. The heat semigroup is strongly continuous on
$\BUC$. The two integral terms in $\mathcal T U$ also depend
continuously on $t$ in this space. Indeed, the part with $0<t-s<\delta$
has norm at most $CK(\sqrt\delta+\delta)$, and the two kernels are
continuous in $L^1$ away from zero time. With
$X_\tau=C([0,\tau];\BUC)$, the kernel estimates yield
\begin{align*}
 \|\mathcal TU-U^0\|_{X_\tau}&\le CK(\sqrt\tau+\tau),\\
 \|\mathcal TU-\mathcal T\widetilde U\|_{X_\tau}
 &\le C\Lambda(\sqrt\tau+\tau)\|U-\widetilde U\|_{X_\tau}.
\end{align*}
Choose $\tau$ so that $CK(\sqrt\tau+\tau)\le r$ and
$C\Lambda(\sqrt\tau+\tau)\le1/2$. Banach's fixed-point theorem gives a
mild solution with $a\ge c$ throughout $[0,\tau]$.
Testing \eqref{eq:mild} against the backward heat operator applied to a
compactly supported smooth function gives \eqref{eq:U-system} in
distributions with initial value $U_{\mathrm{in}}$. The flux and
source are bounded, and $\int_0^\tau\|DG_s\|_1\,ds<\infty$ justifies
interchanging integrals.

For two initial values in a common $\Omega_{A,c}$, take the same
existence interval and common bounds on the two solutions. Subtracting
their mild formulas and using the second estimate above, after shortening
$\tau$ if necessary, gives
\[
 \|U-\widetilde U\|_{X_\tau}
 \le \|U_{\mathrm{in}}-\widetilde U_{\mathrm{in}}\|_\infty
     +\tfrac12\|U-\widetilde U\|_{X_\tau}.
\]
This proves \eqref{eq:stability}. The same estimate applies from any time
at which two positive mild solutions agree and have a common positive
lower bound. By continuity it gives uniqueness on their common interval
in the open set $\{(a,b)\in\BUC:\inf_x a(x)>0\}$. If $U(S)\in\Omega_{A,c}$,
the construction gives an additional interval of length
$\tau_0(A,c)>0$, with the dependence on the fixed parameters understood.
The semigroup law gives \eqref{eq:mild} across the joining time. This
proves the continuation assertion. The regularity assertion follows by
applying Lemma~\ref{lem:smoothing} below to the solution just constructed.
\end{proof}

\subsection{Propagation of spatial regularity}\label{subsec:smoothing}

The local existence time was obtained using only the supremum norm
of $U$ and a positive lower bound for $a$. To apply the energy
identities to mollified data, we need their higher spatial regularity
to persist on this same interval. For an integer $k\ge0$, write
$\BUC^k$ for functions whose spatial
derivatives through order $k$ are bounded and uniformly continuous,
with norm the sum of their supremum norms.

\begin{lemma}\label{lem:smoothing}
Let $U$ be the positive mild solution of \eqref{eq:mild} on $[0,\tau]$
constructed in Appendix~\ref{subsec:proof-local}.
If $U_{\mathrm{in}}\in\BUC^m$, then $U\in C([0,\tau];\BUC^m)$ and,
for $0\le k\le m$,
\begin{equation*}
 \sup_{0\le t\le\tau}\|U(t)\|_{\BUC^k}\le C_k .
\end{equation*}
Here $C_k$ depends on $k$, $\|U_{\mathrm{in}}\|_{\BUC^k}$, $\tau$,
the bounds for $\|U\|_\infty$ and $1/a$ on the interval, and the fixed
parameters. For smooth data with all spatial derivatives bounded, all
space and time derivatives are continuous and bounded on $[0,\tau]\times\R^3$,
with one-sided time derivatives at the endpoints.
\end{lemma}

\begin{proof}
The heat semigroup acts contractively and strongly continuously on
$\BUC^k$, and $\|\partial_jG_t*\|_{\BUC^k\to\BUC^k}\le Ct^{-1/2}$.
For $|\alpha|\le k$,
\begin{equation*}
 D^\alpha\mathcal N[\omega a]=\mathcal N[D^\alpha(\omega a)],\qquad
 \|D^\alpha\mathcal N[\omega a]\|_\infty\le C_k\|a\|_{\BUC^k}.
\end{equation*}
Indeed, the derivatives of $\omega a$ through order $k$ lie in
$L^1\cap\BUC$, by \eqref{eq:weight-powers} and $\omega\in L^1$.
Their translations are continuous in both norms, so taking successive
difference quotients in \eqref{eq:weighted-field} justifies the identity.
The product and chain rules give local Lipschitz continuity of
$\mathcal F_j,\mathcal J$ in $\BUC^k$ on sets with bounded $\BUC^k$
norm and a fixed positive lower bound for $a$. Applying the preceding
contraction argument in $C([0,\delta];\BUC^k)$ produces a local solution
there. Uniqueness in $\BUC$ identifies it with $U$.

Suppose inductively that the $\BUC^{k-1}$ norm is bounded throughout
$[0,\tau]$, where $k\ge1$. Each order-$k$ derivative of a local flux
or source contains derivatives of order $k$ of $U$ only linearly.
The terms of highest order in the field product are
$(D^\alpha a)\mathcal N[\omega a]$ and
$a\mathcal N[\omega D^\alpha a]$. Consequently the differentiated
fluxes and source are bounded in the supremum norm by
$C(1+\|D^kU\|_\infty)$, where $C$ depends only on the induction
bounds, the lower bound for $a$, $k$ and the fixed parameters.
On any subinterval $[t_0,t_0+h]$ of the $\BUC^k$ existence interval,
differentiate the restarted mild formula by transferring the $k$
derivatives to the fluxes, source and initial value. The quantity
$M_h=\sup_{t_0\le s\le t_0+h}\|D^kU(s)\|_\infty$ satisfies
\[
 M_h\le\|D^kU(t_0)\|_\infty+C(\sqrt h+h)(1+M_h).
\]
For a fixed sufficiently small $h_0>0$, this gives
$M_h\le2\|D^kU(t_0)\|_\infty+1$ when $h\le h_0$.
Iteration bounds the $\BUC^k$ norm up to the smaller of $\tau$ and the
maximal $\BUC^k$ existence time. The local contraction can then be
restarted with a uniformly positive lifespan at any earlier time.
This extends the $\BUC^k$ solution through $\tau$ and completes
the induction.

For smooth data the same argument applies for every $k$. The equation
expresses $U_t$ as a continuous $\BUC^k$-valued function for each $k$.
Integrating the distributional equation identifies this function with the
strong time derivative. Repeated time differentiation is justified by
the product and chain rules, and by bounded linearity of
$a\mapsto\mathcal N[\omega a]$ on each $\BUC^k$.
\end{proof}

The initial values $U_{0,n}$ of Corollary~\ref{cor:approx-smooth} have
bounded spatial derivatives of all orders. Lemma~\ref{lem:smoothing} and the positive
lower bound for $a_n$ on each compact time interval in its lifespan
give the claimed regularity for $\rho_n=\omega a_n$, $v_n=b_n/a_n$,
$g_n=d+Da_n/a_n$ and $u_n=v_n-g_n$.

\section*{Acknowledgments}
	X Huang is partially supported by Chinese Academy of Sciences Project for Young Scientists in Basic Research (Grant No. YSBR-031), National Natural Science Foundation of China (Grant Nos. 12494542, 11688101) and National Key R\&D Program of China (Grant No. 2021YFA1000801). 
	
	\vspace{0.3cm}
	\noindent\textbf{AI assistance.} The authors acknowledge the use of AI to explore proof strategies, check calculations, and polish the manuscript. In particular, AI suggested the critical-exponent argument that covers the endpoint case $\gamma=3$ of the polytropic comparison, where the exponent $\gamma-1=2$ of $\|\rho\|_\infty^{\gamma-1}$ matches the density bound $\int_0^T\|\rho(t)\|_\infty^2\,dt\le C_T$, and helped verify the cancellation in the pressure correction. The authors take full responsibility for the mathematical arguments and the final manuscript.

	\vspace{1cm}
	\noindent\textbf{Data availability statement.} Data sharing is not applicable to this article.
	
	\vspace{0.3cm}
	\noindent\textbf{Conflict of interest.} The authors declare that they have no conflict of interest.

\end{document}